\documentclass[a4paper,11pt,leqno]{article}
\usepackage{amsmath,amssymb,latexsym,amscd,amsthm,amsxtra,graphicx}
\begin{document}



\newtheorem{thm}{Theorem}[subsection]
\newtheorem{lem}[thm]{Lemma}
\newtheorem{cor}[thm]{Corollary}
\newtheorem{prop}[thm]{Proposition}
\newtheorem{prob}[thm]{Problem}
\newtheorem{clm}{Claim}
\newtheorem{dfn}[thm]{Definition}


\newcommand{\pf}{{\it Proof} \,\,\,\,\,\,}
\newcommand{\remark}{{\it Remark} \,\,\,\,}

\newcommand{\eqb}{\begin{equation}}
\newcommand{\eqe}{\end{equation}}

\newcommand{\cpx}{{\mathbb C}}
\newcommand{\rel}{{\mathbb R}}
\newcommand{\rat}{{\mathbb Q}}
\newcommand{\itg}{{\mathbb Z}}
\newcommand{\nat}{{\mathbb N}}

\newcommand{\T}{\Theta}
\newcommand{\Th}{\Theta}
\newcommand{\Ta}{\Theta^\alpha}

\newcommand{\www}{{\mathfrak w}}
\newcommand{\gggg}{{\mathfrak g}}
\newcommand{\kkk}{{\mathfrak k}}
\newcommand{\aaa}{{\mathfrak a}}
\newcommand{\aas}{{\mathfrak a}_{\T}}
\newcommand{\aasa}{{\mathfrak a}_{\T^\alpha}}
\newcommand{\bbb}{{\mathfrak b}}
\newcommand{\llll}{{\mathfrak l}}
\newcommand{\llls}{{\mathfrak l}_{\T}}

\newcommand{\ccc}{{\mathfrak c}}
\newcommand{\lls}{{\mathfrak l}_{\T}}
\newcommand{\obb}{\bar{\mathfrak b}}
\newcommand{\aad}{{\mathfrak a}^\ast}
\newcommand{\ads}{{\mathfrak a}^\ast_{\T}}
\newcommand{\nnn}{{\mathfrak n}}
\newcommand{\nns}{{\mathfrak n}_{\T}}
\newcommand{\uuu}{{\mathfrak u}}
\newcommand{\tuu}{\tilde{\mathfrak u}}
\newcommand{\buu}{\bar{\mathfrak u}}
\newcommand{\bns}{\bar{\mathfrak n}_{\T}}
\newcommand{\mmm}{{\mathfrak m}}
\newcommand{\mms}{{\mathfrak m}_{\T}}
\newcommand{\jjj}{{\mathfrak j}}
\newcommand{\hhh}{{\mathfrak h}}
\newcommand{\hhd}{{\mathfrak h}^\ast}
\newcommand{\qqq}{{\mathfrak q}}
\newcommand{\ppp}{{\mathfrak p}}
\newcommand{\pps}{{\mathfrak p}_{\T}}
\newcommand{\sss}{{\mathfrak s}}
\newcommand{\ooo}{{\mathfrak o}}
\newcommand{\Sth}{{\mathfrak S}}
\newcommand{\gl}{{\mathfrak g}{\mathfrak l}}
\newcommand{\so}{\sss\ooo}
\newcommand{\fsp}{{\mathfrak s}{\mathfrak p}}
\newcommand{\hol}{{\cal O}}
\newcommand{\diff}{{\cal D}}
\newcommand{\ana}{{\cal A}}
\newcommand{\cca}{{\cal C}}
\newcommand{\mca}{{\cal M}}
\newcommand{\hca}{{\cal H}}
\newcommand{\gS}{{\mathfrak S}}
\newcommand{\ii}{\sqrt{-1}}
\newcommand{\real}{\mbox{Re}}
\newcommand{\Real}{\mbox{Re}}
\newcommand{\supp}{\mbox{supp}}
\newcommand{\rank}{\mbox{rank}}
\newcommand{\card}{\mbox{card}}
\newcommand{\Ad}{\mbox{Ad}}
\newcommand{\ad}{\mbox{ad}}
\newcommand{\dg}{{\deg}}
\newcommand{\Hom}{\mbox{Hom}}
\newcommand{\End}{\mbox{End}}
\newcommand{\Ext}{\mbox{Ext}}
\newcommand{\fgt}{\mbox{Fgt}}
\newcommand{\Dim}{\mbox{Dim}}
\newcommand{\gr}{\mbox{gr}}
\newcommand{\Ann}{\mbox{Ann}}
\newcommand{\tr}{\mbox{tr}}
\newcommand{\hht}{\mbox{ht}}
\newcommand{\Mod}{{\sf Mod}}
\newcommand{\sgn}{\mbox{sgn}}
\newcommand{\pro}{\mbox{pro}}
\newcommand{\Ind}{\mbox{\sf Ind}}
\newcommand{\SO}{\mbox{SO}}
\newcommand{\Oo}{\mbox{O}}
\newcommand{\SOo}{{\mbox{SO}_0}}
\newcommand{\GL}{\mbox{GL}}
\newcommand{\Spp}{\mbox{Sp}}
\newcommand{\U}{\mbox{U}}
\newcommand{\triv}{{\mbox{triv}}}
\newcommand{\Rea}{{\mbox{Re}}}

\newcommand{\PP}{{\sf P}}
\newcommand{\PS}{{\sf P}_{\T}^{++}}
\newcommand{\PSP}{{\sf P}_{\T^\prime}^{++}}
\newcommand{\PPS}{{\sf P}_{\T}^{++}}
\newcommand{\PPSA}{{\sf P}_{\T^\alpha}^{++}}
\newcommand{\PPSK}{{\sf P}_{\T^k}^{++}}
\newcommand{\Q}{{\sf Q}}

\newcommand{\lel}{\stackrel{L}{\leq}}
\newcommand{\lell}{\stackrel{L}{\leq}}
\newcommand{\ler}{\stackrel{R}{\leq}}
\newcommand{\lerr}{\stackrel{R}{\leq}}
\newcommand{\lelr}{\stackrel{LR}{\leq}}
\newcommand{\gel}{\stackrel{L}{\geq}}
\newcommand{\ger}{\stackrel{R}{\leq}}
\newcommand{\gelr}{\stackrel{LR}{\leq}}
\newcommand{\lleq}{\stackrel{L}{\sim}}
\newcommand{\rreq}{\stackrel{R}{\sim}}
\newcommand{\lreq}{\stackrel{LR}{\sim}}
\newcommand{\leqs}{\leqslant}
\newcommand{\geqs}{\geqslant}


\title{{\sf The homomorphisms between scalar generalized Verma modules
associated to maximal parabolic subalgebras}}
\author{{\sf Hisayosi Matumoto}
\\Graduate School of Mathematical Sciences\\ University of Tokyo\\ 3-8-1
Komaba, Tokyo\\ 153-8914, JAPAN\\ e-mail: hisayosi@ms.u-tokyo.ac.jp}
\date{}
\maketitle
\begin{abstract}
Let $\gggg$ be a finite-dimensional simple Lie algebra over $\cpx$.
We classify the homomorphisms between $\gggg$-modules induced  from one-dimensional
 modules of maximal parabolic subalgebras.\footnote{Keywords: 
 generalized Verma module, differential invariant
 \\ AMS Mathematical Subject Classification: 22E47, 17B10}
\end{abstract}

\setcounter{section}{0}
\section*{\S\,\, 0.\,\,\,\, Introduction}
\setcounter{subsection}{0}

In this article, we consider the existence problem of homomorphisms between generalized
Verma modules, which are induced from one dimensional representations
(such generalized Verma modules are called scalar, cf.\   [Boe 1985]).
Our main result is the classification of the  homomorphisms between
scalar generalized Verma modules with respect to the  maximal parabolic subalgebras.

A sufficient condition for the existence of the homomorphisms between
Verma modules is given by [Verma 1968].
Bernstein, I.\ M.\ Gelfand, and S.\ I.\ Gelfand proved the condition of
Verma is also a necessary condition. ([Bernstein-Gelfand-Gelfand 1975]) 

Later, Lepowsky studied the problem for the generalized Verma modules.
In particular, Lepowsky ([Lepowsky 1975a]) solved the existence problem
of nontrivial homomorphisms between scalar generalized Verma modules
associated to the parabolic
subalgebras which are the complexifications of the minimal parabolic
subalgebras of real rank one simple Lie algebras (so-called the real rank one
case).

Lepowsky also obtained a sufficient condition for the
existence of the homomorphisms between scalar generalized Verma modules
associated to the complexification of the minimal parabolic subalgebras
of (not necessarily rank one) real semisimple Lie algebras.
His condition is quite similar to that of Verma and he conjectured it is
also a sufficient condition in the setting of  complexified minimal
parabolic algebras ([Lepowsky 1975b]).

Boe ([Boe 1985]) solved the existence problem in the case  of parabolic subalgebras whose nilradical is
commutative (so-called the Hermitian symmetric case). 

The existence problem for maximal parabolic algebras is, in principle,
reduced to the Kazdhan-Lusztig algorithm. 
Casian and Collingwood ([Casian-Collingwood 1987]) proposed a direct
method of computing the Kazdhan-Lusztig data involving the generalized
Verma modules.
Applying the Kazdhan-Lusztig algorithm works very well in some cases.
The structures of the (not necessarily scalar) generalized
Verma modules in the real rank one case and in the Hermitian symmetric
case  are studied precisely ([Boe-Collingwood 1985],
[Boe-Enright-Shelton 1988], [Collingwood-Irving-Shelton 1988]).
Boe and Collingwood ([Boe-Collingwood 1990]) studied the case that all
the (not necessarily scalar) generalized Verma modules are multiplicity
free.
The cases treated by Boe and Collingwood is more general than the real rank
one case and the  Hermitian symmetric case (so-called the multiplicity
free case).
In that case, they studied the structures of   the (not necessarily scalar) generalized
Verma modules precisely.
In particular, they solved the existence problem  of the non-trivial
homomorphisms for regular integral
infinitesimal characters  in the multiplicity
free case.

However, I surmise it is not easy to give the explicit answer to the existence problem by the
Kazdhan-Lusztig algorithm in the
general setting. 
Our central dogma is ``{\it Consider the most singular parameter, then everything
turns to be easy.}''
Our approach to the problem consists of the following three main
ingredients.

(1) \,\,\, The translation principle

The translation principle has a long history. 
In [Vogan 1988], Vogan proposed  an idea on translation principle in
order to establish the irreducibility of a discrete series representation
of a semisimple symmetric space in some case.
His idea is extremely useful for the study of the existence problem.
 Depending Vogan's idea, we formulated a version of translation principle in
[Matumoto 1993] Proposition 2.2.3.
In some cases, this enable us to reduce  the existence of a non-trivial
homomorphism to the most singular case in which
 the problem is often trivial. 

(2) \,\,\, Jantzen's irreducibility criterion

Applying a version of the translation principle,  we can often
reduce  the nonexistence of  a non-trivial
homomorphism to the irreducibility of a  particular generalized Verma
module.
In [Jantzen 1977], Jantzen gave a sufficient and necessary condition for
the irreducibility of a generated Verma module.
His result is extremely useful for our purpose and we can establish the
nonexistence of nontrivial homomorphisms in many cases.

(3) \,\,\, The Kazdhan-Lusztig theory ([Kazdhan-Lusztig 1979],
    [Brylinski-Kashiwara 1981], [Beilinson-Bernstein 1993])

Although we do not compute Kazdhan-Lusztig polynomials, the existence of
the Kazdhan-Lusztig algorithm plays an important role in our approach.
The point is that the definition of the Kazdhan-Lusztig polynomials only
depends on Coxeter systems.
In some cases, this enable us to reduce the problem to that of a
different maximal parabolic subalgebra  of a different simple Lie
algebra, which is easier than the original problem.

In the most of the cases we can solve the existence problem by the above three ideas.
However, in some cases we need extra arguments.

This article consists of five sections.

We fix notations and introduce some fundamental material in \S 1.

In \S 2, we introduce  sufficient conditions for the existence problem.

In \S 3, we treat the case of the classical algebras.
The type A case is in the Hermitian symmetric case. So, we only
consider the case of the type B,C, and D.
The main theorems are Theorem 3.2.1, 3.2.2, and 3.2.3.

In \S 4, we treated the case of the exceptional algebras.

In \S 5, using a comparison result, we explain how to construct a
homomorphism between scalar generalized Verma modules associated to  a general
parabolic subalgebra  from a homomorphism between scalar generalized
Verma modules associated to a maximal
parabolic subalgebra.
We call such a homomorphism an elementary homomorphism.

I would like to propose:

{\bf Working Hypothesis} \,\,\, {\it An arbitrary nontrivial homomorphism between  scalar
generalized Verma modules is a composition of elementary homomorphisms. }

The working hypothesis in the case of the Verma modules is nothing but the result
of Bernstein-Gelfand-Gelfand.
The first statement of the Lepowsky conjecture ([Lepowsky 1975b]
Conjecture 6.13) means that, in the case of complexified
minimal parabolic subalgebras of real semisimple Lie algebras, the above
working hypothesis is affirmative. 

The result of this article solves the existence of edge-of-wedge type
embeddings in the case of  the maximal parabolic subgroups of complex
reductive groups
 ([Matumoto 2003]).

\subsection*{Acknowledgements}

During the proceeding of this study, I participated in the
``Representation Theory of Lie Groups''  at the Institute for Mathematical
Sciences (IMS) the National University of Singapore.
It is my pleasure to thank their warm hospitality.


\setcounter{section}{1}
\setcounter{subsection}{0}

\section*{\S\,\, 1.\,\,\,\,Notations and Preliminaries }

\subsection{General notations}

In this article, we use the following notations and conventions.

As usual we denote the complex number field, the real number field, the
ring of (rational) integers, and the set of non-negative integers by
$\cpx$, $\rel$, $\itg$, and $\nat$ respectively.
$\frac{1}{2}\nat$ means the set $\left\{\left. \frac{n}{2}\right|
n\in\nat\right\}$, and
$\frac{1}{2}+\nat$ means the set $\left\{\left. \frac{1}{2}+n\right| n\in\nat\right\}$.
We denote by $\emptyset$ the empty set.
For any (non-commutative) $\cpx$-algebra $R$, ``ideal'' means  ``2-sided ideal'', ``$R$-module'' means ``left $R$-module'', and sometimes we denote by $0$ (resp.\ $1$) the trivial $R$-module $\{0\}$ (resp.\ $\cpx$).
Often, we identify a (small) category and the set of its objects.
Hereafter ``$\dim$'' means the dimension as a complex vector space, and ``$\otimes$'' (resp. $\Hom$) means the tensor product over $\cpx$ (resp. the space of $\cpx$-linear mappings), unless we specify.
For a complex vector space $V$, we denote by $V^\ast$ the dual vector space.
For $a,b\in\cpx$, ``$a\leqslant b$'' means that $a,b\in\rel$ and $a\leqslant b$.
We denote by $A- B$ the set theoretical difference.
$\card A$ means the cardinality of a set $A$.

\subsection{Notations for reductive Lie algebras}

Let $\gggg$ be a complex reductive Lie algebra, $U(\gggg)$ the universal
enveloping algebra of $\gggg$, and  $\hhh$ a Cartan subalgebra of
$\gggg$.
We denote by $\Delta$ the root system with respect to
$(\gggg,\hhh)$.
We fix some positive root system $\Delta^+$ and let $\Pi$ be the set of simple roots.
Let $W$ be the Weyl group of the pair $(\gggg, \hhh)$ and let
$\langle\,\,,\,\,\rangle$ be a non-degenerate invariant bilinear form on
$\gggg$.
For $w\in W$, we denote by $\ell(w)$ the length of $w$ as usuall.
We also denote the inner product on $\hhd$ which is induced from the
above form by the same symbols $\langle\,\,,\,\,\rangle$.
For $\alpha\in\Delta$, we denote by $s_\alpha$ the reflection in $W$ with respect to $\alpha$.
We denote by $w_0$ the longest element of $W$.
For $\alpha\in\Delta$, we define the coroot $\check{\alpha}$ by
$\check{\alpha} = \frac{2\alpha}{\langle\alpha,\alpha\rangle}$,
as usual.
We call $\lambda\in \hhd$ is dominant (resp. anti-dominant), if $\langle\lambda,\check{\alpha}\rangle$ is not a negative (resp. positive) integer, for each $\alpha\in \Delta^+$.
We call $\lambda\in\hhd$ regular, if $\langle\lambda,\alpha\rangle\neq 0$, for each $\alpha\in\Delta$.
We denote by $\PP$ the integral weight lattice, namely
$\PP = \{\lambda\in\hhd\mid\langle\lambda,\check{\alpha}\rangle\in\itg \,\,\,\,\hbox{for all}\,\,\, \alpha \in\Delta\}$.
If $\lambda\in\hhd$ is contained in $\PP$, we call $\lambda$ an integral weight.
We define $\rho\in\PP$ by $\rho = \frac{1}{2}\sum_{\alpha\in\Delta^+}\alpha$.
Put 
$\gggg_\alpha =\{X\in\gggg\mid\forall H\in\hhh \;\; [H,X] =\alpha(H)X\}$,
$\uuu = \sum_{\alpha\in\Delta^+}\gggg_\alpha$,
$\bbb = \hhh+\uuu$.
Then $\bbb$ is a Borel subalgebra of $\gggg$.
We denote by $\Q$ the root lattice, namely $\itg$-linear span of $\Delta$.  We also denote by $\Q^+$ the linear combination of $\Pi$ with non-negative integral coefficients.  
 For $\lambda\in\hhh^\ast$, we denote by $W_\lambda$ the integral Weyl group.
Namely,
\[ W_\lambda=\{w\in W\mid w\lambda-\lambda\in\Q\}.\]
We denote by $\Delta_\lambda$ the set of integral roots.
\[\Delta_\lambda=\{\alpha\in\Delta\mid \langle\lambda,\check{\alpha}\rangle\in\itg\}.\]
It is well-known that $W_\lambda$ is the Weyl group for $\Delta_\lambda$.
We put $\Delta^+_\lambda=\Delta^+\cap\Delta_\lambda$.
This is a positive system of $\Delta_\lambda$.  We denote by $\Pi_\lambda$ the set of simple roots for $\Delta^+_\lambda$ and denote by $\Phi_\lambda$ the set of reflection corresponding to the elements in $\Pi_\lambda$.  
So, $(W_\lambda, \Phi_\lambda)$ is a Coxeter system.
We denote by $\Q_\lambda$ the integral root lattice, namely $\Q_\lambda=\itg\Delta^+_\lambda$ and put $\Q^+_\lambda=\nat\Pi_\lambda$.

Next, we fix notations for a parabolic subalgebra (which contains $\bbb$).
Hereafter, through this article we fix an arbitrary subset $\T$ of $\Pi$.
Let $\bar{\T}$ be the set of the elements of $\Delta$ which are written by linear combinations of elements of $\T$ over $\itg$.
Put
$\aas = \{H\in\hhh\mid \forall \alpha\in \T \,\,\,\alpha(H)=0\}$,
$\lls = \hhh +\sum_{\alpha\in\bar{\T}}\gggg_\alpha$,  
$\nns = \sum_{\alpha\in\Delta^+\backslash\bar{\T}}\gggg_\alpha$,  
$\pps = \lls +\nns$.
Then $\pps$ is a parabolic subalgebra of $\gggg$ which contains $\bbb$.
Conversely, for an arbitrary parabolic subalgebra $\ppp\supseteq\bbb$, there exists some $\T \subseteq \Pi$ such that $\ppp =\pps$.
We denote by $W_\T$ the Weyl group for $(\lls,\hhh)$. 
$W_\T$ is identified with a subgroup of $W$ generated by $\{s_\alpha\mid \alpha\in\T\}$.
We denote by $w_\T$ the longest element of $W_\T$.
Using the invariant non-degenerate bilinear form $\langle\,\,,\,\,\rangle$, we regard  ${\aaa_{\T}}^\ast$ as a subspace of $\hhh^\ast$. 
It is known that there is a unique nilpotent (adjoint) orbit (say $\hol_{\ppp_\T}$) whose intersection with $\nns$ is Zarisky dense in $\nns$.
$\hol_{\pps}$ is called the Richardson orbit with respect to $\pps$.
We denote by $\bar{\hol}_{\ppp_{\T}}$ the closure of $\hol_{\ppp_{\T}}$ in $\gggg$.
Put
$\rho_{\T} = \frac{1}{2}(\rho-w_{\T}\rho)$ and 
$\rho^{\T} = \frac{1}{2}(\rho+w_{\T}\rho)$.
Then, $\rho^{\T}\in{\aaa_{\T}}^\ast$.

\subsection{Generalized Verma modules}

Define
\begin{align*}
\PS & = \{\lambda\in\hhd\mid \forall\alpha\in \T\,\,\,\,\,
 \langle\lambda,\check{\alpha}\rangle\in\{1,2,...\}\}\\
{}^\circ\PS & =\{\lambda\in\hhd\mid \forall\alpha\in \T\,\,\,\,\,
 \langle\lambda,\check{\alpha}\rangle=1\}
\end{align*}
We easily have
\begin{align*}
{}^\circ\PS =\{\rho_{\T}+\mu\mid \mu\in\ads\}.
\end{align*}
For $\mu\in\hhd$ such that $\mu+\rho\in\PS$, we denote by $\sigma_\T(\mu)$ the irreducible finite-dimensional $\lls$-representation whose highest weight is $\mu$.
Let $E_\T(\mu)$ be the representation space of $\sigma_\T(\mu)$.
We define a left action of $\nns$ on $E_\T(\mu)$ by $X\cdot v =0$ for all $X\in\nns$ and $v\in E_\T(\mu)$.
So, we regard $E_\T(\mu)$ as a $U(\pps)$-module.

For $\mu\in\PS$, we define a generalized Verma module ([Lepowsky 1977]) as follows.
\[ M_\T(\mu) = U(\gggg)\otimes_{U(\pps)}E_\T({\mu-\rho}).\]
For all $\lambda\in\hhd$, we write $M(\lambda) = M_\emptyset(\lambda)$.
$M(\lambda)$ is called a Verma module.
For $\mu\in\PS$, $M_\T(\mu)$ is a quotient  module of $M(\mu)$.
Let $L(\mu)$ be the unique highest weight $U(\gggg)$-module with the highest weight $\mu-\rho$.
Namely, $L(\mu)$ is a unique irreducible quotient of $M(\mu)$.
For $\mu\in\PS$, the canonical projection of $M(\mu)$ to
$L(\mu)$ is factored by $M_\T(\mu)$.

$\dim E_{\T}(\mu-\rho)=1$ if and only if $\mu\in {}^\circ\PS$.
If $\mu\in {}^\circ\PS$, we call $M_\T(\mu)$ a scalar generalized
Verma module.

For a finitely generated $U(\gggg)$-module $V$, we denote by $\Dim(V)$
(resp.\ $c(V)$) the Gelfand-Krillov dimesion (resp.\ the multiplicity)
of $V$. (See [Vogan 1978]).
We easily see $\Dim(M_{\T}(\mu))=\dim\nns$ and
$c((M_{\T}(\mu))=\dim E_\T({\mu-\rho})$.

The following result is one of the fundamental results on the existence
problem of homomorphisms between scalar generalized Verma modules.

\begin{thm} ([Lepowsky 1976])

Let $\mu,\nu\in {}^\circ\PS$.

(1) \,\,\,\, $\dim \Hom_{U(\gggg)}(M_{\T}(\mu), M_{\T}(\nu))\leqslant 1.$

(2) \,\,\,\, Any non-zero homomoorphism of $M_{\T}(\mu)$ to
 $M_{\T}(\nu)$ is injective.
\end{thm}

Hence, the existence problem of  homomorphisms  between scalar
generalized Verma modules is reduce to the following problem.

{\bf Problem} \,\, Let $\mu,\nu\in {}^\circ\PS$.
 When is $M_{\T}(\mu)\subseteq M_{\T}(\nu)$ ?

\subsection{Homomorphisms associated with Duflo involutions}

Herafter we assume  $\T\subseteq\Pi$, $\mu\in \PS$ and $\mu$ is dominant and regular.
Then, we easily have $w_{\T}w_0\mu\in \PS$. and $M_{\T}(w_{\T}w_0\mu)$ is
irreducible.
Here, we consider the following problem

\begin{prob}

When is $M_{\T}(w_{\T}w_0\mu)\hookrightarrow M_{\T}(\mu)$?
\end{prob}

Concerning to Problem 1.4.1, a necessary and sufficient condition is known.

\begin{thm} \,\,\,\, ([Matumoto 1993])
Let $\T\subseteq\Pi$ and $\mu\in \PS$.
If $\mu$ is dominant and regular, then 
 the following two conditions (1) and (2) are equivalent.

(1) \,\,\,\, $M_{\T}(w_{\T}w_0\mu)\hookrightarrow M_{\T}(\mu)$.

(2) \,\,\,\, $w_{\T}w_0$ is a Dulfo involution for the Coxeter system
 $(W_{\mu}, \Phi_{\mu})$.

\end{thm}

In particular, the answer of Problem 2.1.1 only
depend on the Coxter system $(W_{\mu}, \Phi_{\mu})$.
(In fact, this fact is a conclusion of  the Kazdhan-Lusztig conjecture.)
We can find a complex reductive Lie algebra whose Weyl group (with a
set of the simple reflections) are isomorphic to $(W_{\mu}, \Phi_{\mu})$.
So, we can deduce Problem 1.4.1 to the
following special case.

\begin{prob}

When is $M_{\T}(w_{\T}w_0\rho)\hookrightarrow M_{\T}(\rho)$?
\end{prob}

We also remark that the following easy fact.
\begin{lem}
Let $\T\subseteq\Pi$ and $\lambda\in {}^\circ\PPS$ such that
 $w_{\T}w_0=w_0w_{\T}$  and $w_{\T}w_0\in W_\lambda$.
We denote by $w^\prime_0$ be the longest element of $W_\lambda$ with
 respect to $\Pi_\lambda$.
Then, $w_0=w_0^\prime$ and $w_{\T}w_0^\prime=w^\prime_0 w_{\T}$.
\end{lem}

\subsection{Translation principle and its application }

We denote by $Z(\gggg)$ the center of $U(\gggg)$.
It is well-known that $Z(\gggg)$ acts on $M(\lambda)$ by the Harish-Chandra homomorphism $\chi_\lambda : Z(\gggg) \rightarrow \cpx$ for all $\lambda$.
$\chi_\lambda =\chi_\mu$ if and only if there exists some $w\in W$ such that $\lambda =w\mu$.
We denote by $\hbox{\bf Z}_\lambda$ the kernel of $\chi_\lambda$ in $Z(\gggg)$.
Let $M$ be a $U(\gggg)$-module and $\lambda\in\hhh^\ast$.
We say that $M$ has an infinitesimal character $\lambda$ if and only if $Z(\gggg)$ acts on $M$ by $\chi_\lambda$.
We say that $M$ has a generalized infinitesimal character $\lambda$ if and only if for any $v\in M$ there is some positive integer $n$ such that ${\hbox{\bf Z}_\lambda}^nv=0$.
  We say $M$ is locally $Z(\gggg)$-finite, if and only if for any $v\in M$ we have  $\dim Z(\gggg)v<\infty$.
We denote by $\mca_{Zf}$ (cf.\ [Bernstein-Gelfand 1980]) the category of $Z(\gggg)$-finite $U(\gggg)$-modules.
We also denote by $\mca[\lambda]$ the category of $U(\gggg)$-modules with generalized infinitesimal character $\lambda$.
Then, from the Chinese remainder theorem, we have a direct sum of abelian categories 
$\mca_{Zf}=\bigoplus_{\lambda\in\hhh^\ast}\mca[\lambda]$.
We denote by $P_\lambda$ the projection functor from $\mca_{Zf}$ to $\mca[\lambda]$.
For $\mu\in\PP$, we denote by $V_\mu$ the irreducible finite-dimensional $U(\gggg)$-module with an extreme weight $\mu$.
Let $\mu,\lambda\in\hhh^\ast$ satisfy $\mu-\lambda\in\PP$.  
Let $M$ be an object of $\mca[\lambda]$.
Then, from a result of Kostant we have that $M\otimes V_{\mu-\lambda}$ is an object of $\mca_{Zf}$.  
So, we can define translation functor $T^\mu_\lambda$ from $\mca[\lambda]$ to $\mca[\mu]$ as follows.
\[T^\mu_\lambda(M)=P_\mu(M\otimes V_{\mu-\lambda}).\]
\bigskip
The translation functors are exact.

We put 
\[ W({\T})=\{w\in W\mid w\T= \T \}.\]
Then, $W({\T})$ is a subgroup of $W$.  
Moreover, $w\rho_{\T}=\rho_{\T}$ and $w_{\T}w=ww_{\T}$ hold for all $w\in W({\T})$ and $W({\T})$ preserves $\aas^\ast$.
In particular, $W({\T})\subseteq W_{\rho_{\T}}$.

We say that $\lambda\in\ads$ is strongly $\T$-antidominant if and only if $\langle\lambda,\alpha\rangle\leqslant 0$  for all $\alpha\in\Q^+\cap
\Q_{\rho_S+\lambda}$. (Cf.\ [Matumoto 1993])

Next, we consider the images of generalized Verma modules under certain translation functors.

\begin{lem} (Cf.\ [Matumoto 1993] Lemma 1.2.3)
Assume that $\mu,\lambda\in\ads$ are strongly ${\T}$-antidominant and that $\lambda-\mu$ is dominant and integral.
Let $w\in W({\T})$.
Then, we have
\[ T^{-\rho_{\T}+\lambda}_{-\rho_{\T}+\mu}(M_{\T}(\rho_{\T}+w\mu))=M_{\T}(\rho_{\T}+w\lambda).\]

\end{lem}


\setcounter{section}{2}
\setcounter{subsection}{0}

\section*{\S\,\, 2.\,\,\,\,Sufficient conditions }

For almost all $(\gggg,\hhh, \Delta^+, \T)$, the necessary and
sufficient condition given in Theorem 1.4.2 is hard to check.
So, we consider sufficient conditions, which we can check easily.

\subsection{A sufficient condition}

We fix the notations for characters. (Cf.\ [Dixmier 1977] 7.5.1, [Knapp
2002] V. 6)
Let $\cpx^{\hhh^\ast}$ be the $\cpx$-vector space of all functions from
$\hhh^\ast$ to $\cpx$.
For $f\in\cpx^{\hhh^\ast}$, we define $\supp(f)=\{\lambda\in\hhd\mid
f(\lambda)\neq 0\}.$
For $\lambda\in\hhd$, we define $e^\lambda$ the member of
$\cpx^{\hhh^\ast}$, that is $1$ at $\lambda$ and $0$ elsewhere.
Let $\cpx<\hhd>$ be the set of all $f\in\cpx^{\hhh^\ast}$ such that
$\supp(f)$ is contained in the union of a finite number of sets
$\nu_i-Q^+$ with each $\nu_i$ in $\hhd$.
We introduce the structure of $\cpx$-algebra on $\cpx^{\hhh^\ast}$ as in [Knapp 2002] (5.65).

Let $V$ be a $U(\gggg)$-module. 
For $\lambda\in\hhd$, we define the weight space with respect to
$\lambda$ as follows.
\begin{align*}
V_\lambda=\{v\in V\mid \forall H\in\hhh \,\,\, Hv=\lambda(H)v \}
\end{align*}

We say that $V$ has a character if $V$ is the direct sum of its weight
spaces under $\hhh$ and if $\dim V_\lambda\leqslant \infty$ for all
$\lambda\in\hhd$.
In this case, the character is
\begin{align*}
[V]=\sum_{\lambda\in\hhd}(\dim V_\lambda)e^\lambda.
\end{align*}

For example, for $\lambda\in\PPS$, the following formula is well-known.
\begin{align*}
[M_{\T}(\lambda)] = D^{-1}\sum_{w\in W_\T}(-1)^{\ell(w)}e^{w\lambda}.
\end{align*}
Here, we denote by $D$ the Weyl denominator, namely
$D=e^\rho\prod_{\alpha\in\Delta^+}(1-e^{-\alpha})$.
In particular, we have
$[M(\lambda)] = D^{-1}e^{\lambda}$.

We put ${\bf B}_{st}=\{[M(w\rho)]\mid w\in W\}$ and denote by $\cca$
the subspace of $\cpx^{\hhh^\ast}$ spanned by ${\bf B}_{st}$.
We also put ${\bf B}_{irr}=\{[L(w\rho)]\mid w\in W\}$.
Then,  ${\bf B}_{st}$ and ${\bf B}_{irr}$ are bases of $\cca$.
We identify the group algebra $\cpx[W]$ with $\cca$ via the
correspondence $W\ni w\leftrightsquigarrow [M(ww_0\rho)]\in{\bf
B}_{st}$.
We can introduce a $W$-module structure on $\cca$ identifying $\cca$ and
the right regular representation on $\cpx[W]$.
We consider the equivalence relations $\lleq$ and $\rreq$ on $W$ defined
in [Kazdhan-Lusztig 1979].
The representation theoretic meanings of these equivalence relations are as follows.

\begin{thm} ([Joseph 1977], [Vogan 1980])
Let $x,y\in W$. Then we have

(1) \,\,\,\, $x\lleq y$ if and only if
 $\Ann_{U(\gggg)}(L(xw_0\rho))=\Ann_{U(\gggg)}(L(yw_0\rho))$.

(2) \,\,\,\, $x\rreq y$ if and only if
there exist some finite-dimensional $U(\gggg)$-modules $E_1$
 and $E_2$ such that $L(xw_0\rho)$ and $L(yw_0\rho)$ are irreducible
 constituents of $E_1\otimes L(yw_0\rho)$ and $E_2\otimes L(xw_0\rho)$,
 respectively.
\end{thm}

For $x\in W$, we denote by $V_x^R$ the $\cpx$-vector space with a basis
$B_x=\{[L(yw_0\rho)]\mid y\rreq x\}$.
If we identify $V_x^R$ with a subquotient of $\cca$ appropriately, we
may regard  $V_x^R$ a $W$-module. (For example, see [Barbasch-Vogan
1983]).
$V^R_x$ is called a right cell representation.

We denote by $\hca$ the space of $W$-harmonic polynomials on $\hhd$,
which can be regarded as a $W$-module in a usual manner (cf.\ [Vogan
1978]).

We quote:

\begin{thm} ([Vogan 1978], [Joseph 1980a 1980b])

For $x\in W$,
there is a $W$-homomorphism $\phi_x$ of $V^R_x$ to $\hca$ satisfying the
 following conditions.

(1) \,\,\,\, The image $\phi_x(V^R_x)\subseteq \hca$ is an irreducible
 representation. 
 
(2) \,\,\,\,  The image $\phi_x(V^R_x)$ is the special representation
 corresponding to the unique open dense nilpotent orbit in the
 associated variety of $\Ann_{U(\gggg)}(L(xw_0\rho))$ via the Springer correspondence.
In particular,  $\phi_{w_\T}(V^R_{w_\T})$ is the special representation
 corresponding to the Richardson orbit $\hol_{\ppp_\T}$.

(3) \,\,\,\, For $y,z\in W$ such that $y\rreq x\rreq z$, $y\lleq z$ if
 and only if $\phi_x([L(yw_0\rho)])$  and $\phi_x([L(zw_0\rho)])$ are
 proportional to each other.

\end{thm}

Now, we state a sufficient condition. 

\begin{prop}
Assume that $\T\subseteq\Pi$ satisfy $w_{\T}w_0=w_0w_{\T}$.
Moreover, we assume that $V^R_{w_{\T}}$ is irreducible as a $W$-module. 
Then, we have $M_{\T}(w_{\T}w_0\rho)\subseteq M_{\T}(\rho)$.
\end{prop}

\proof

We put $I=\Ann_{U(\gggg)}(M_{\T}(w_{\T}w_0\rho))$.
$I$ is primitive, since $M_{\T}(w_{\T}w_0\rho)$ is irreducible.
Since $w_{\T}w_0=w_0w_{\T}$, we have $w_0w_{\T}\in W(S)$.
So, from [Borho-Jantzen 1977] 4.10 Corollary, we have
$I=\Ann_{U(\gggg)}(M_{\T}(\rho))$.
Since $c(M_{\T}(\rho))=1$, $M_{\T}(\rho)$ has a unique irreducible
constituent of maximal Gelfand-Killirov dimension (say $L(\sigma
\rho)$). Here, $\sigma$ is an element of $W$.
(In fact we we easy to prove $L(\sigma
\rho)$ is the unique irreducible submodule of $M_{\T}(\rho)$.  )
Assuming that $\sigma\neq w_{\T}w_0$, we shall deduce a contradiction.

When we regard  a $U(\gggg)$-module $E$ as a $U(\pps)$-module, we write it by $E|_{\pps}$.
Since $M_{\T}(\mu)\otimes E\cong
U(\gggg)\otimes_{U(\pps)}(E_{\T}(\mu-\rho)\otimes E|_{\pps})$ holds, we
easily see that $\sigma w_0\rreq w_{\T}$.

From [Borho-Kraft 1976] 3.6, we see $I=\Ann_{U(\gggg)}(L(\sigma\rho))$.
From Theorem 1.6.1, we have $w_{\T}\lleq \sigma w_0$.

Hence, from Theorem 1.6.2, $\phi_{w_{\T}}$ is not injective.
This contradicts our assumption that $V^R_{w_{\T}}$ is  irreducible.
\,\,\,\, Q.E.D.

The multiplicity of a special representation $\phi_x(V^R_x)$ in
a right cell $V^R_x$ is one.
Moreover, any irreducible constituents in the right cell $V^R_x$
belongs to the same family (see [Lusztig 1984] p78) as $\phi_x(V^R_x)$.
So, we have:

\begin{cor}
Assume that $\T\subseteq\Pi$ satisfy $w_{\T}w_0=w_0w_{\T}$ and that
 the family of the special representation corresponding to
 the Richardson orbit $\hol_{\ppp_{\T}}$ does not contain any other
 element.
Then, we have  $M_{\T}(w_{\T}w_0\rho)\subseteq M_{\T}(\rho)$.
\end{cor}

We denote by $\lreq$ the equivalence relation on $W$ generated by
$\lleq$ and $\rreq$.
The following result is well-known and follows from Theorem 2.1.1.
\begin{cor} Let $x,y\in W$ be such that $x\lreq y$.
Then $\Dim(L(xw_0\rho))=\Dim(L(yw_0\rho))$.
\end{cor}
From [Barbasch-Vogan 1983] Corollary 2.24 implies that $x\lreq y$ if
and only if $xw_0\lreq yw_0$.
Hence, we have:
\begin{lem}
Let $x,y\in W$ be such that $x\lreq y$.
Then $\Dim(L(x\rho))=\Dim(L(y\rho))$.
\end{lem}

\subsection{Maximal parabolic subalgebras}

Hereafter we fix $\alpha\in\Pi$.
Put $\T^\alpha=\Pi-\{\alpha \}$.
$\aasa$ is one-dimensional and spanned by $\rho^{\Ta}$.
Moreover, we have
\begin{align*}
{}^\circ\PPSA=\{\rho_{\Ta}+t\rho^{\Ta}\mid t\in\cpx\}.
\end{align*}

We denote by $\omega_\alpha$ the fundamental weight corresponding to $\alpha$.
For any $\beta\in\Ta=\Pi-\{\alpha\}$, we have
$\langle\beta,\rho^{\Ta}\rangle=0$. 
Hence there exists some $d_\alpha\in\rel$ such that
$d_\alpha\omega_\alpha=\rho^{\Ta}$.
Since $2\rho^{\Ta}$ is integral, we have $d_\alpha\in \frac{1}{2}\nat$.

Fir simplicity, for $t\in\cpx$,  we write $M_{\Ta}[t]$ for $M_{\Ta}(\rho_{\Ta}+t\omega_\alpha)$.
We have:

\begin{lem}
 Let $s$ and $t$ be distinct complex numbers such that
 $M_{\Ta}[s]\subseteq M_{\Ta}[t]$.
Then, $s=-t$.
\end{lem}

\proof
Since $M_{\Ta}[s]$ and $M_{\Ta}[t]$ have the same infinitesimal
character, there exists some $w\in W$ such that
$\rho_{\Ta}+t\omega_\alpha=w(\rho_{\Ta}+s\omega_\alpha)$.
Hence, $\langle\rho_{\Ta}+t\omega_\alpha,\rho_{\Ta}+t\omega_\alpha\rangle=
\langle\rho_{\Ta}+s\omega_\alpha,\rho_{\Ta}+s\omega_\alpha\rangle$.
From $\langle\omega_\alpha,\rho_{\Ta}\rangle=0$, we have $t^2=s^2$.
So, $s=-t$.\,\,\,\, $\Box$

We easily have:

\begin{lem}

(1) If $w_{\Ta}w_0=w_0w_{\Ta}$, then
 $w_{\Ta}w_0(\rho_{\Ta}+t\omega_\alpha)=\rho_{\Ta}-t\omega_\alpha$ for all
 $t\in\cpx$.

(2) If $w_{\Ta}w_0\neq w_0w_{\Ta}$, then
 $w_{\Ta}w_0(\rho_{\Ta}+t\omega_\alpha)=\rho_{\Ta}-t\omega_\alpha$ if and only
 if $t=d_\alpha$.
\end{lem}

We put 
\begin{align*}
c_\alpha=\min\{c\in\rel \mid 2c \omega_\alpha\in \Q^+\}. 
\end{align*}

Clearly $2c_\alpha$ is a positive integer.

\begin{lem}
If $w_{\Ta}w_0=w_0w_{\Ta}$, then either $c_\alpha=1$ or
 $c_\alpha=\frac{1}{2}$.
\end{lem}

\proof
We have only to show $c_\alpha\leqslant 1$.
We may assume $\gggg$ is simple.
If $\gggg$ is a simple Lie algebra of the type other than $A_n$,
$D_{2n+1}$, and $E_6$, then the exponent of $\Q/\PP$ is $1$ or $2$.
So, in this case $2\omega_\alpha\in \Q^+$ for any $\alpha\in\Pi$.
For the case of the type $A_n$,
$D_{2n+1}$, or $E_6$, we can check $2\omega_\alpha\in \Q^+$ under the
assumption $w_{\Ta}w_0=w_0w_{\Ta}$ by the case-by-case analysis.
\,\,\,\, $\Box$

If  $M_{\Ta}[-t]\subseteq M_{\Ta}[t]$ holds, then
$2t\omega_\alpha=(\rho_{\Ta}+t\omega_\alpha)-(\rho_{\Ta}-t\omega_\alpha)$
is in $\Q^+$.
So, we have:
\begin{cor}
If $w_{\Ta}w_0=w_0w_{\Ta}$ and 
if  $M_{\Ta}[-t]\subseteq M_{\Ta}[t]$ holds, then
 $t\in\frac{1}{2}\nat$.
Moreover, if $c_\alpha=1$, then $M_{\Ta}[-t]\hookrightarrow M_{\Ta}[t]$
 implies 
 $t\in\nat$.
\end{cor}
\begin{dfn}
If $t\in\frac{1}{2}\nat$ and $\rho_{\Ta}+t\omega_\alpha$ is not
integral, we say $\rho_{\Ta}+t\omega_\alpha$ is half-integral.
\end{dfn}

We examine behavior of the translation functors in the setting of this subsection. 

First, Lemma 1.5.1 and the exactness of the translation functor imply:
\begin{lem}
Assume $w_{\Ta}w_0=w_0w_{\Ta}$.
Let  $t\in\rel$ and $n\in\nat$  be such that  $t-n\geqslant 0$.
Then, $M_{\Ta}[-t]\subseteq M_{\Ta}[t]$ implies
 $M_{\Ta}[-t+n]\subseteq M_{\Ta}[t-n]$.
\end{lem}

From the translation principle, we also have:

\begin{lem}
Assume $w_{\Ta}w_0=w_0w_{\Ta}$.
Let  $t\in\rel$ be such that  $\rho_{\Ta}+t\omega_\alpha$ is dominant regular.
Let $n\in\nat$.
Then, $M_{\Ta}[-t]\subseteq M_{\Ta}[t]$ implies
 $M_{\Ta}[-t-n]\subseteq M_{\Ta}[t+n]$.
\end{lem}

In case $\rho_{\Ta}+t\omega_\alpha$ is not dominant regular, the corresponding
statement to Lemma 1.7.6 is not necessarily correct.
In fact, we need an extra assumption.

Let $G$ be a complex connected reductive Lie group, whose Lie algebra is
$\gggg$.
Let $P_{\T}$ be the parabolic subgroup of $G$ corresponding to
$\ppp_{\T}$.
We consider the generalized flag variety $X_{\T}=G/P_{\T}$.
Since the holomorphic cotangent bundle $T^\ast X_\T$ has a natural
symplectic structure, we can construct the moment map $m_{\T}:
T^\ast X_{\T}\rightarrow \gggg^\ast$.
Using $\langle\, , \,\rangle$, we identify $\gggg$ and $\gggg^\ast$.
Then, we regard the moment map as a surjective map of $T^\ast X_{\T}$ to
the closure of the Richardson orbit $\hol_{\ppp_{\T}}$.

We easily see the number $N_1(P_{\T})$ defined in [Hesselink 1978]
1.4 Step3 is the degree of the moment map $m_{\T} : T^\ast
X_{\T}\rightarrow  \overline{\hol_{\ppp_{\T}}}$.
So, $N_1(P_{\T})=1$ if and only if $m_{\T}$ is birational.
For classical Lie algebras, $N_1(P_{\T})$ is obtained in [Hesselink
1978] 7.1 Theorem.

\begin{lem} ([Matumoto 1993] Proposition 2.2.3)
\mbox{}
Assume $w_{\Ta}w_0=w_0w_{\Ta}$ and the moment map $m_{\T} : T^\ast X_{\Ta}\rightarrow \overline{\hol_{\ppp_{\Ta}}}$ is birational.
Let $n\in\nat$.
Then, $M_{\Ta}[-t]\subseteq M_{\Ta}[t]$ implies
 $M_{\Ta}[-t-n]\subseteq M_{\Ta}[t+n]$.
\end{lem}

Since $M_{\Ta}[0]\subseteq M_{\Ta}[0]$, we have
\begin{cor}([Matumoto 1993] Corollary 2.2.4)
\mbox{}

Assume $w_{\Ta}w_0=w_0w_{\Ta}$ and the moment map $m_{\T} : T^\ast X_{\Ta}\rightarrow \overline{\hol_{\ppp_{\Ta}}}$ is birational.
Then, we have $M_{\Ta}[-n]\subseteq M_{\Ta}[n]$ for all $n\in\nat$.
\end{cor}

Next, we introduce Jantzen's criterion for the irreducibility of a
generalized Verma module.

For any $\lambda\in\hhd$, we define an element of $\cpx<\hhd>$ as follows.
\[
 \Upsilon_{\Ta}(\lambda)=D^{-1}\sum_{w\in W_{\Ta}}(-1)^{\ell(w)}e^{w\lambda}.
\]
Of course, for $\lambda\in\PPSA$, we have
$\Upsilon_{\Ta(\lambda)}=[M_{\Ta}(\lambda)]$.

Put $\Delta_{\Ta}=\itg\Ta\cap\Delta$ and
$\left(\Delta^{\Ta}\right)^+=\Delta^+-\Delta_{\Ta}$.
We immediately have:
\begin{cor} \mbox{}

(1) \,\, If $\lambda \in\hhd$ satisfies $\langle\beta,\lambda\rangle=0$
 for some $\beta\in\Delta_{\Ta}$,
 then $\Upsilon_{\Ta}(\lambda)=0$.

(2) \,\, For $\lambda \in\hhd$ and $w\in W_{\Ta}$, we have
$\Upsilon_{\Ta}(\lambda)-(-1)^{\ell(w)}\Upsilon_{\Ta}(w\lambda)=0$.
\end{cor}

The following result is a special case of [Jantzen 1977] Satz 3.

\begin{thm} \, (Jantzen)
\mbox{}

For $\lambda\in\PPSA$, the following (1) and (2) are equivalent.

(1) \,\,\,\, $M_{\Ta}(\lambda)$ is irreducible.

(2) \,\,\,\, We have
\begin{align*}
\sum_{\substack{\beta\in\left(\Delta^{\Ta}\right)^+ \\
 \langle\lambda,\beta^\vee\rangle\in\nat-\{0\}
}}\Upsilon_{\Ta}(s_\beta\lambda)=0.
\end{align*}
\end{thm}
 
{\bf Remark} \,\,\, The above statement is slightly different from that in
Jantzen's paper.
However, we consider maximal parabolic subalgebras.
In this case, we can easily see the above condition (2) is equivalent to
the condition of Jantzen.


\setcounter{section}{3}
\setcounter{subsection}{0}

\section*{\S\,\, 3.\,\,\,\, Classical Lie algebras}

Throughout this section, $n$ means a positive integer such that
$n\geqslant 2$ (resp.\ $n\geqslant 4$) whenever we consider the simple
Lie algebra of the type
$B_n$ or $C_n$ (resp.\ $D_n$).

\subsection{The root systems}

We retain the notations in \S 1 and \S 2.

{\bf (${\bf B_n}$ type)} \,\,\,\,
We consider the root system $\Delta$ for $\gggg=\so(2n+1,\cpx)$.
Then we can choose an orthonormal basis $e_1,...,e_n$ of $\hhh^\ast$ such that 
\begin{align*}
\Delta = \{ \pm e_i\pm e_j\mid 1\leqslant i< j\leqslant
n\}\cup
\{\pm e_i\mid 1\leqslant i\leqslant n\}.
\end{align*}
We choose a positive system as follows.
\begin{align*}
\Delta^+ = \{ e_i\pm e_j\mid 1\leqslant i< j\leqslant
n\}\cup
\{e_i\mid 1\leqslant i\leqslant n\}.
\end{align*}
If we put $\alpha_i=e_i-e_{i+1}$  \,\,\, $(1\leqslant i <n)$ and
$\alpha_n=e_n$,
then $\Pi=\{\alpha_1,...,\alpha_n\}$.

For simplicity, we write $\T^k$ for $\T^{\alpha_k}$
 for $1\leqslant k\leqslant n$.
Then $\llll_{\T^k}$ is isomorphic to $\gl(k,\cpx)\times
\so(2(n-k)+1,\cpx)$.
Since $w_0$ is contained in the center of $W$, we have
 $w_0w_{\T^k}=w_{\T^k}w_0$ for any $1\leqslant k\leqslant n$.

We write $\omega_k$, $c_k$, and $d_k$ for $\omega_{\alpha_k}$,
 $c_{\alpha_k}$, and $d_{\alpha_k}$, respectively.
Then, for $1\leqslant k<n$, we have $\omega_k=e_1+\cdots+e_k$, $c_k=\frac{1}{2}$, and
 $d_k=n-\frac{k}{2}$.
If $k=n$, then $\omega_n=\frac{1}{2}(e_1+\cdots+e_n)$, $c_n=1$, and $d_n=n$.

Assume $k$ is even or $k=n$.
Then,  $\rho_{\T^k}+t\omega_k$ is integral if and only if
 $t\in\itg$.

Assume $k$ is odd and $1\leqslant k<n$.
Then, $\rho_{\T^k}+t\omega_k$ is integral if and only if
 $t-\frac{1}{2}\in\itg$.

{\bf (${\bf C_n}$ type)} \,\,\,\,

We consider the root system $\Delta$ for $\gggg=\sss\ppp(n,\cpx)$.
Then we can choose an orthonormal basis $e_1,...,e_n$ of $\hhh^\ast$ such that 
\begin{align*}
\Delta = \{ \pm e_i\pm e_j\mid 1\leqslant i< j\leqslant
n\}\cup
\{\pm 2e_i\mid 1\leqslant i\leqslant n\}.
\end{align*}
We choose a positive system as follows.
\begin{align*}
\Delta^+ = \{ e_i\pm e_j\mid 1\leqslant i< j\leqslant
n\}\cup
\{2e_i\mid 1\leqslant i\leqslant n\}.
\end{align*}
If we put $\alpha_i=e_i-e_{i+1}$  \,\,\, $(1\leqslant i <n)$ and
$\alpha_n=2e_n$,
then $\Pi=\{\alpha_1,...,\alpha_n\}$.

For simplicity, we write $\T^k$ for $\T^{\alpha_k}$
 for $1\leqslant k\leqslant n$.
Then $\llll_{\T^k}$ is isomorphic to $\gl(k,\cpx)\times
\sss\ppp(n-k,\cpx)$.
Since $w_0$ is contained in the center of $W$, we have
 $w_0w_{\T^k}=w_{\T^k}w_0$ for any $1\leqslant k\leqslant n$.

We write $\omega_k$, $c_k$, and $d_k$ for $\omega_{\alpha_k}$,
 $c_{\alpha_k}$, and $d_{\alpha_k}$, respectively.
Then, for $1\leqslant k\leqslant n$ we have $\omega_k=e_1+\cdots+e_k$ and
 $d_k=n-\frac{k-1}{2}$.
If $k$ is even (resp.\ odd), then $c_k=\frac{1}{2}$ (resp.\ $c_k=1$).

Assume $k$ is odd.
Then,  $\rho_{\T^k}+t\omega_k$ is integral if and only if
 $t\in\itg$.

Assume $k$ is even.
Then, $\rho_{\T^k}+t\omega_k$ is integral if and only if
 $t-\frac{1}{2}\in\itg$.

{\bf (${\bf D_n}$ type)} \,\,\,\,

We consider the root system $\Delta$ for $\gggg=\so(2n,\cpx)$.
Then we can choose an orthonormal basis $e_1,...,e_n$ of $\hhh^\ast$ such that 
\begin{align*}
\Delta = \{ \pm e_i\pm e_j\mid 1\leqslant i< j\leqslant
n\}.
\end{align*}
We choose a positive system as follows.
\begin{align*}
\Delta^+ = \{ e_i\pm e_j\mid 1\leqslant i< j\leqslant
n\}.
\end{align*}
If we put $\alpha_i=e_i-e_{i+1}$  \,\,\, $(1\leqslant i <n)$ and
$\alpha_n=e_{n-1}+e_n$,
then $\Pi=\{\alpha_1,...,\alpha_n\}$.

For simplicity, we write $\T^k$ for $\T^{\alpha_k}$
 for $1\leqslant k\leqslant n$.
Since the case of $k=n-1$ is essentially same as the case of $k=n$,
 since $\ppp_{\T^{n-1}}$ and $\ppp_{\T^{n}}$ are conjugate under an automorphism of $\gggg$.
So, when we consider the type D case, we omitt the case of $k=n-1$.

Then, $\llll_{\T^k}$ is isomorphic to $\gl(k,\cpx)\times
\so(2(n-k),\cpx)$.
We have
 $w_0w_{\T^k}=w_{\T^k}w_0$ for any $1\leqslant k <n$.
However, $w_0w_{\T^n}=w_{\T^n}w_0$ if and only if $n$ is even.

We write $\omega_k$, $c_k$, and $d_k$ for $\omega_{\alpha_k}$,
 $c_{\alpha_k}$, and $d_{\alpha_k}$, respectively.
Then, for $1\leqslant k<n-1$, we have $\omega_k=e_1+\cdots+e_k$ and
 $d_k=n-\frac{k+1}{2}$.
If $k$ is even (resp.\ odd) and $k\neq n$, then $c_k=\frac{1}{2}$ (resp.\ $c_k=1$).

For $k=n$, we have $\omega_n=\frac{1}{2}(e_1+\cdots+e_n)$, $c_n=1$, and $d_n=n$.
If $n$ is even (resp.\ odd), then $c_n=1$ (resp.\ $c_k=2$).

Assume $k$ is odd or $k=n$.
Then,  $\rho_{\T^k}+t\omega_k$ is integral if and only if
 $t\in\itg$.

Assume $k$ is even and $1\leqslant k<n$.
Then, $\rho_{\T^k}+t\omega_k$ is integral if and only if
 $t-\frac{1}{2}\in\itg$.

\subsection{Statements of the main results}

Here, we describe the existence of homomorphisms between scalar
generalized Verma modules with respect to the maximal parabolic
subalgebras of classical Lie algebras.
For simple Lie algebras of type A, the answer is given in [Boe 1985].
So, we treat the case of the types of B,C, and D. 

\begin{thm} ($B_n$-type)

Consider the case of $\gggg=\so(2n+1,\cpx)$.
Let $k$ be a positive integer such that $k\leqslant n$.

(1) \,\,\, We consider the case that $3k< 2n+1$.
Then, $M_{\T^k}[-t]\subseteq M_{\T^k}[t]$ if and only if $t\in\nat$. 

(2) \,\,\, We consider the case that $3k\geqslant 2n+1$.

\mbox{} \,\, (2a) \,\, Assume $k$ is odd and $k\neq n$. Then, 
 $M_{\T^k}[-t]\subseteq M_{\T^k}[t]$ if and only if
 $t\in\frac{1}{2}\nat$.

\mbox{} \,\, (2b) \,\, Assume $n$ is odd. Then, 
 $M_{\T^n}[-t]\subseteq M_{\T^n}[t]$ if and only if
 $t\in\nat$.

\mbox{} \,\, (2c) \,\, Assume $k$ is even. Then,
 $M_{\T^k}[-t]\subseteq M_{\T^k}[t]$ if and only if $t=0$.

\end{thm}

{\bf Remark} \,\,\,\, The case of $k=1$ is due to [Lepowsky 1975a].
The cases of $k=1,2,3,n$ are in the multiplicity free case in
[Boe-Collingwood 1990].

\begin{thm} ($C_n$-type)

Consider the case of $\gggg=\sss\ppp(n,\cpx)$.
Let $k$ be a positive integer such that $k\leqslant n$.

(1) \,\, We consider the case that $3k\leqslant 2n$.

\mbox{} \,\, (1a) \,\, Assume $k$ is odd.  
Then, $M_{\T^k}[-t]\subseteq M_{\T^k}[t]$ if and only if $t=0$.
 
\mbox{} \,\, (1b) \,\, Assume $k$ is even.
Then, $M_{\T^k}[-t]\subseteq M_{\T^k}[t]$ if and only if $t\in\frac{1}{2}\nat$.

(2) \,\, We consider the case that  $3k>2n$. Then,
 $M_{\T^k}[-t]\subseteq M_{\T^k}[t]$ if and only if
 $t\in\nat$.

\end{thm}

{\bf Remark} \,\,\,\, The case of $k=2$ is due to [Lepowsky 1975a].
The case of $k=n$ is due to [Boe 1985].
The cases of $k=1,2,3,n$ are in the multiplicity free case in
[Boe-Collingwood 1990]. 

\begin{thm} ($D_n$-type)

Consider the case of $\gggg=\so(2n,\cpx)$.
Let $k$ be a positive integer such that $k\leqslant n-2$ or $k=n$.

(1) \,\,\, We consider the case that $3k< 2n$.
Then, $M_{\T^k}[-t]\subseteq M_{\T^k}[t]$ if and only if $t\in\nat$. 

(2) \,\,\, We consider the case that  $3k\geqslant 2n$.
 
\mbox{} \,\, (2a) \,\, Assume $k$ is odd. Then, 
 $M_{\T^k}[-t]\subseteq M_{\T^k}[t]$ if and only if
 $t=0$.

\mbox{} \,\, (2b) \,\, Assume $k$ is even and $k\neq n$. Then,
 $M_{\T^k}[-t]\subseteq M_{\T^k}[t]$ if and only if $t\in\frac{1}{2}\nat$.

\mbox{} \,\, (2c) \,\, Assume $n$ is even. Then,
 $M_{\T^n}[-t]\subseteq M_{\T^n}[t]$ if and only if $t\in\nat$.

\end{thm}

{\bf Remark} \,\,\,\, The case of $k=1$ is due to [Lepowsky 1975a].
The case of $k=n$ is due to [Boe 1985].
The cases of $k=1,2,n$ are in the multiplicity free case in
[Boe-Collingwood 1990]. 

\medskip

We give proofs of the above theorems in the subsquent sections.

\subsection{Richardson orbits}

We consider a partition $\pi=(p_1,...,p_k)$ of a positive integer $m$
 such that $0<p_1\leqslant p_2\leqslant \cdots\leqslant
p_k$ and $p_1+p_2+\cdots+p_k=m$.
We put $\pi[i]=\card\{j\mid p_j=i\}$ for any positive integer $i$.
Let $\{h_1,...,h_r\}=\{i\in\nat\mid i>0, \pi[i]>0 \}$. 
Then, we write $h_1^{\pi[h_1]}\cdots h_r^{\pi[h_r]}$ for $\pi$.
For example, $5\cdot 3^2\cdot 1^5 $ means $(1,1,1,1,1,3,3,5)$. 

{\bf Type ${\bf B_n}$}
\,\,\,\,
The nilpotent orbits in $\gggg=\so(2n+1,\cpx)$ are parametrized by
partitions $\pi$ of $2n+1$ such that, for any
even number $2i$,  $\pi[2i]$ is even. (For example see
[Carter 1985] p394)
From [Collingwood-McGovern 1993] 7 and [Hesselink 1978] 7.1, we
have:
\begin{lem} Let $\gggg=\so(2n+1,\cpx)$.

(1) \,\,\,\, If $3k< 2n+1$, then the Richardson orbit $\hol_{\ppp_{\T^k}}$
 corresponds to the partition $3^k\cdot 1^{2n+1-3k}$.
In this case, the moment map $m_{\T^k}$ is birational.

(2) \,\,\,\, If $3k\geqslant 2n+1$ and $k$ is odd, then  the Richardson orbit 
$\hol_{\ppp_{\T^k}}$
 corresponds to the partition $3^{2n+1-2k}\cdot 2^{3k-2n-1}$.
In this case, the moment map $m_{\T^k}$ is birational.

(3) \,\,\,\, If $3k\geqslant 2n+1$ and $k$ is even, then  the Richardson orbit 
$\hol_{\ppp_{\T^k}}$
 corresponds to the partition $3^{2n+1-2k}\cdot 2^{3k-2n-2}\cdot 1^2$.
\end{lem}

From Corollary 2.2.9, and Lemma 3.3.1, we
have:
\begin{lem}  Let $\gggg=\so(2n+1,\cpx)$.
Assume that $3k\leqslant 2n+1$ or $k$ is odd.
Then, we have $M_{\T^k}[-t]\subseteq M_{\T^k}[t]$ for all $t\in\nat$.
\end{lem}

{\bf Type ${\bf C_n}$}
\,\,\,\,
The nilpotent orbits in $\gggg=\sss\ppp(n,\cpx)$ are parametrized by
partitions $\pi$ of $2n$ such that, for any
odd number $2i+1$,  $\pi[2i+1]$ is even. (For example see
[Carter 1985] p394)
From [Collingwood-McGovern 1993] 7 and [Carter 1985] Chapter 13, we
have:
\begin{lem} Let $\gggg=\sss\ppp(2n,\cpx)$.

(1) \,\,\,\, If $3k\leqslant 2n$ and $k$ is even, then the Richardson orbit $\hol_{\ppp_{\T^k}}$
 corresponds to the partition $3^k\cdot 1^{2n-3k}$.
In this case the moment map $m_{\T^k}$ is birational.

(2) \,\,\,\, If $3k\leqslant 2n$ and $k$ is odd, then the Richardson orbit $\hol_{\ppp_{\T^k}}$
 corresponds to the partition $3^{k-1}\cdot 2^2\cdot 1^{2n-3k-1}$.

(3) \,\,\,\, If $3k>2n$, then  the Richardson orbit 
$\hol_{\ppp_{\T^k}}$
 corresponds to the partition $3^{2n-2k}\cdot 2^{3k-2n}$.
In this case, the moment map $m_{\T^k}$ is birational.
\end{lem}

From  Corollary 2.2.9, and Lemma 3.3.3, we
have:
\begin{lem}  Let $\gggg=\sss\ppp(2n,\cpx)$.
Assume that $3k>2n$  or $k$ is even.
Then, we have $M_{\T^k}[-t]\subseteq M_{\T^k}[t]$ for all $t\in\nat$.
\end{lem}

{\bf Type ${\bf D_n}$}
\,\,\,\,
We can associate a nilpotent orbit in $\gggg=\so(2n,\cpx)$ with 
a partition $\pi=(p_1,...,p_k)$ of $2n$ such that, for any
even number $2i$,  $\pi[2i]$ is even.
If there is at least one odd number in  $p_1,...,p_k$, there is a unique
nilpotent orbit associated with  $(p_1,...,p_k)$.  
The exceptional orbits are so-called very even nilpotent orbits.
If $p_1,...,p_k$ are all even, there are two nilpotent orbits
associated with the partition $(p_1,...,p_k)$.
(For example see
[Carter 1985] p395)
From [Collingwood-McGovern 1993] 7 and [Hesselink 1978] 7.1, we
have:
\begin{lem} Let $\gggg=\so(2n,\cpx)$ and let $1\leqslant k\leqslant n-2$
 or $k=n$.

(1) \,\,\,\, If $3k\leqslant 2n$, then the Richardson orbit $\hol_{\ppp_{\T^k}}$
 corresponds to the partition $3^k\cdot 1^{2n-3k}$.
In this case,  the moment map $m_{\T^k}$ is birational.

(2) \,\,\,\, If $3k>2n$ and $k$ is odd, then  the Richardson orbit 
$\hol_{\ppp_{\T^k}}$
 corresponds to the partition $3^{2n-2k}\cdot 2^{3k-2n-1}\cdot 1^2$.

(3) \,\,\,\, If $3k>2n$ and $k$ is even, then  the Richardson orbit 
$\hol_{\ppp_{\T^k}}$
 corresponds to the partition $3^{2n-2k}\cdot 2^{3k-2n}$.
In this case,  the moment map $m_{\T^k}$ is birational.
\end{lem}

From Corollary 2.2.9, and Lemma 3.3.5, we
have:
\begin{lem}  Let $\gggg=\so(2n,\cpx)$.
Assume that $3k\leqslant 2n$ or $k$ is even.
Then, we have $M_{\T^k}[-t]\subseteq M_{\T^k}[t]$ for all $t\in\nat$.
\end{lem}

\subsection{Existence results via comparison}

In this subsection, we prove the following result.

\begin{lem}\,\,\, \mbox{}

(1) \,\,\,\, Assume that $\gggg=\sss\ppp(n,\cpx)$ and that $k$ is an
 even positive integer such that $3k\leqslant 2n$.
Then, we have $M_{\T^k}[-t]\subseteq M_{\T^k}[t]$ for all $t\in\frac{1}{2}+\nat$.

(2) \,\,\,\, Assume that $\gggg=\so(2n+1,\cpx)$ and that $k$ is an
odd positive integer such that $3k\geqslant 2n+1$ and $k\neq n$.
Then, we have $M_{\T^k}[-t]\subseteq M_{\T^k}[t]$ for all $t\in\frac{1}{2}+\nat$.
\end{lem}

{\it Proof} \,\,\,\,
We prove (1). 
We can prove (2) in a similar way.

Let $k$ be an
 even positive integer such that $3k\leqslant 2n$.
At first, we consider the case that $\gggg=\so(2n+1,\cpx)$.
In this case $\rho^{\T^k}=d_k\omega_k=(n-\frac{k}{2})\omega_k$ from 3.1.
Therefore,  $\rho=\rho_{\T^k}+(n-\frac{k}{2})\omega_k$ and
$w_{\T^k}w_0\rho=\rho_{\T^k}-(n-\frac{k}{2})\omega_k$.
Since $n-\frac{k}{2}\in\nat$, Lemma 3.3.2 implies
$M_{\T^k}(w_{\T^k}w_0\rho)\subseteq M_{\T^k}(\rho)$.
From Theorem 1.4.2, $w_{\T^k}w_0$ is a Duflo involution in the Weyl
group for $(\gggg,\hhh)$, where $\gggg=\so(2n+1,\cpx)$.
However, the Weyl group of the $B_n$-type and that of the $C_n$-type
 are isomorphic to each other as a  Coxeter system.
Since the notion of the Duflo involutions only depends on the structure
 of the Coxeter system, we see that as an element of the Weyl group for
 $\sss\ppp(n, \cpx)$, $w_{\T^k}w_0$ is a Duflo involution.
From Theorem 1.4.2 implies $M_{\T^k}(w_{\T^k}w_0\rho)\subseteq
 M_{\T^k}(\rho)$ holds for the case that $\gggg=\sss\ppp(n,\cpx)$.
Namely, $M_{\T^k}[-d_k]\subseteq M_{\T^k}[d_k]$.
In this case, $d_k=n-\frac{n-1}{2}\in \frac{1}{2}+\nat$.
So, from Lemma 2.2.6 and Lemma 2.2.7, we have $M_{\T^k}[-t]\subseteq M_{\T^k}[t]$ for all $t\in\frac{1}{2}+\nat$. \,\,\,\, Q.E.D.

\subsection{Irreducibility of a right cell}

In this subsection, we treat the remaining case that nontrivial homomorphisms exist.
Namely, we show:
\begin{lem}
 \,\,\,\, Assume that $\gggg=\so(2n,\cpx)$ and that $k$ is an
even positive integer such that $3k\geqslant 2n$, $k\leqslant n$, and $k\neq n-1$.
Then, we have $M_{\T^k}[-t]\subseteq M_{\T^k}[t]$ for all $t\in\frac{1}{2}+\nat$.
\end{lem}

\proof
If $k=n$, then $\nnn_{\T^k}$ is abelian.
This case is treated in [Boe 1985].
So, we assume that $k$ is an
even positive integer such that $3k>2n$, $k<n-1$.
Put $s=\frac{k}{2}$.
As in [Carter 1985] p376, we can associate an irreducible representation
of the Weyl group of the type $D_n$ with a so-called symbol.
A symbol is a pair of strictly increasing sequence of non-negative integers of the
same length.
We identify two symbols;
\[
\left( 
\begin{array}{r}
\lambda_1,...,\lambda_k\\
\mu_1,...,\mu_k
\end{array}
\right) 
\,\,\,\,\, \mbox{and} \,\,\,\,\,
\left( 
\begin{array}{r}
\mu_1,...,\mu_k \\
\lambda_1,...,\lambda_k
\end{array}
\right). 
\]
As in [Carter 1985] 13.3, we can associate the partition $3^{2n-2k}\cdot
2^{3k-2n}$ to a pair of partitions as follows.
\[
\left( 
\begin{array}{r}
1,2,3,..., 3s-n,3s-n+2,3s-n+3,3s-n+4,...,s+2 \\
1,2,3,..., 3s-n,3s-n+1,3s-n+2,3s-n+3,...,s+1
\end{array}
\right) 
\]
This is the symbol associated with the special representation (say $\pi_k$)
corresponding to the Richardson orbit $\hol_{\ppp_{\T^k}}$ via the
Springer correspondence. (Lusztig, Shoji)

As in [Carter 1985] 13.2,
two irreducible representation with symbols:
\[
\left( 
\begin{array}{r}
\lambda_1,...,\lambda_k\\
\mu_1,...,\mu_k
\end{array}
\right) 
\,\,\,\,\, \mbox{and} \,\,\,\,\,
\left( 
\begin{array}{r}
\lambda_1^\prime,...,\lambda_k^\prime\\
\mu_1^\prime,...,\mu_k^\prime
\end{array}
\right) 
\]
are contained in the same family if and only if 
$\{\lambda_1,...,\lambda_k,\mu_1,...,\mu_k\}=\{\lambda_1^\prime,...,\lambda_k^\prime,\mu_1^\prime,...,\mu_k^\prime\}$.
So, we easily see the family of $\pi_k$  consists of only one element.
From Corollary 2.1.4, we have $M_{\T^k}(w_{\T^k}w_0\rho)\subseteq
M_{\T^k}(\rho)$.
In this case, $d_k=n-\frac{k+1}{2}\in\frac{1}{2}+\nat$.
So, from Lemma 2.2.6 and Lemma 2.2.7, we have 
$M_{\T^k}[-t]\subseteq M_{\T^k}[t]$ for all $t\in\frac{1}{2}+\nat$.
\,\,\,\, Q.E.D.

{\it Remark}  \,\,\,\, We can prove Lemma 3.4.1 in a similar way to
Lemma 3.4.2.

\subsection{Irreducibility of some generalized Verma modules}

We have:
\begin{lem} \mbox{}

(1) \,\,\,\, Assume that $\gggg=\so(2n+1,\cpx)$ and that $k$ is an even
 positive integer such that $k<n$.
Then, $M_{\T^k}(\frac{1}{2})$ is
 irreducible.

(2) \,\,\,\, Assume that $\gggg=\so(2n+1,\cpx)$ and that $n$ is even.
Then, $M_{\T^n}(1)$ is
 irreducible.
 
(3) \,\,\,\, Assume that $\gggg=\so(2n,\cpx)$ and that $k$ is an even
 positive integer such that $3k<2n$.
Then, $M_{\T^k}(\frac{1}{2})$ is
 irreducible.

(4) \,\,\,\, Assume that $\gggg=\so(2n+1,\cpx)$ and that $k$ is an odd
 positive integer such that $3k<2n+1$.
Then, $M_{\T^k}(\frac{1}{2})$ is
 irreducible.

(5) \,\,\,\, Assume that $\gggg=\sss\ppp(n,\cpx)$ and that $k$ is an even
 positive integer such that $3k>2n$.
Then, $M_{\T^k}(\frac{1}{2})$ is
 irreducible.
 
(6) \,\,\,\, Assume that $\gggg=\sss\ppp(n,\cpx)$.
Then, $M_{\T^1}(1)$ is
 irreducible.
\end{lem}

\proof
For $2\leqslant r\leqslant k$, we put $r^\star=k+2-r$.

First, we prove (1). 
We put $s=\frac{k}{2}$.  
Then, we see
\begin{align*}
\rho_{\T^k}+\frac{1}{2}\omega_k=\sum_{i=1}^{2s}(s-i+1)e_i +
 \sum_{j=2s+1}^n (n-j+\frac{1}{2})e_j.
\end{align*}

We easily see
\begin{multline*}
\left\{\beta\in \left(\Delta^{\Ta}\right)^+\left|
 \langle\rho_{\T^k}+\frac{1}{2}\omega_k,\beta^\vee\rangle\in\nat-\{0\}\right.\right\}\\
=\{e_i\mid
 1\leqslant i \leqslant s\}\cup\{e_i+e_j\mid 1\leqslant i<j\leqslant 2s,
 2s+2>i+j\}.
\end{multline*}

For $2\leqslant i \leqslant s$, $\langle e_i-e_{2s-i+2},
s_{e_i}(\rho_{\T^k}+\frac{1}{2}\omega_k)\rangle=0$.
Hence, Corollary 2.2.10 implies
$\Upsilon_{\Ta}(s_{e_i}(\rho_{\T^k}+\frac{1}{2}\omega_k))=0$.

Next, we assume that $1\leqslant i<j\leqslant 2s$ and that 
 $2s+2>i+j$.
Since $i<j$, $2s+2>i+j$ implies $i\leqslant s$. 
If $2\leqslant i$, then we see $j\neq i^\star$ and $\langle e_j-e_{i^\star},
 s_{e_i+e_j}(\rho_{\T^k}+\frac{1}{2}\omega_k)\rangle =0$.
So,  Corollary 2.2.10 implies $\Upsilon_{\Ta}(s_{e_i+e_j}(\rho_{\T^k}+\frac{1}{2}\omega_k))=0$.
If $i=1$ and $j\neq s+1$, then $j\neq j^\star\neq 1$ and $\langle e_1-e_{j^\star},
 s_{e_i+e_j}(\rho_{\T^k}+\frac{1}{2}\omega_k)\rangle =0$.
So,  Corollary 2.2.10 implies $\Upsilon_{\T^k}(s_{e_i+e_j}(\rho_{\T^k}+\frac{1}{2}\omega_k))=0$.

Hence, we have
\begin{align*}
\sum_{\substack{\beta\in\left(\Delta^{\T^k}\right)^+ \\
 \langle\rho_{\T^k}+\frac{1}{2}\omega_k,\beta^\vee\rangle\in\nat-\{0\}
}}\Upsilon_{\T^k}(s_\beta(\rho_{\T^k}+\frac{1}{2}\omega_k))=\Upsilon_{\T^k}(s_{e_1}(\rho_{\T^k}+\frac{1}{2}\omega_k))+\Upsilon_{\T^k}(s_{e_1+e_{s+1}}(\rho_{\T^k}+\frac{1}{2}\omega_k)).
\end{align*}

Since $\langle e_{s+1},\rho_{\T^k}+\frac{1}{2}\omega_k\rangle=0$, we
have
$s_{e_1}(\rho_{\T^k}+\frac{1}{2}\omega_k)=s_{e_1-e_{s+1}}s_{e_1+e_{s+1}}(\rho_{\T^k}+\frac{1}{2}\omega_k)$.
Hence from Corollary 2.2.10, we have 
\begin{align*}
\sum_{\substack{\beta\in\left(\Delta^{\T^k}\right)^+ \\
 \langle\rho_{\T^k}+\frac{1}{2}\omega_k,\beta^\vee\rangle\in\nat-\{0\}
}}\Upsilon_{\T^k}(s_\beta(\rho_{\T^k}+\frac{1}{2}\omega_k))=0.
\end{align*}
From Jantzen's criterion (Theorem 2.2.11), $M_{\T^k}[\frac{1}{2}]$ is
irreducible.
So, we get (1).

(2) is proved in the same way as (1).

Next, we prove (3).
 
We put $s=\frac{k}{2}$.  
Then, we see
\begin{align*}
\rho_{\T^k}+\frac{1}{2}\omega_k=\sum_{i=1}^{2s}(s-i+1)e_i +
 \sum_{j=2s+1}^n (n-j)e_j.
\end{align*}
We easily see
\begin{multline*}
\left\{\beta\in \left(\Delta^{\T^k}\right)^+\left|
 \langle\rho_{\T^k}+\frac{1}{2}\omega_k,\beta^\vee\rangle\in\nat-\{0\}\right.\right\}\\
=\{e_i+e_j\mid 1\leqslant i<j\leqslant 2s,
 2s+2>i+j\}\\
\cup \{e_i\pm e_j\mid 1\leqslant i\leqslant 2s<j\leqslant n, (s-i+1)\pm(n-j)>0\}.
\end{multline*}
Since $3k<2n$, we have $2s<n-s$.
Let $\beta\in \left(\Delta^{\T^k}\right)^+$ be such that 
 $\langle\rho_{\T^k}+\frac{1}{2}\omega_k,\beta^\vee\rangle\in\nat-\{0\}$.
Moreover, we assume $\beta$ is neither $e_1+e_{s+1}$ nor $e_1+e_{n-s}$.
Then, we can easily see there exists some $\gamma\in\Delta_{\T^k}$ such
that $\langle
s_\beta(\rho_{\T^k}+\frac{1}{2}\omega_k),\gamma^\vee\rangle=0$.
From Corollary 2.2.10, we have
\begin{align*}
\sum_{\substack{\beta\in\left(\Delta^{\T^k}\right)^+ \\
 \langle\rho_{\T^k}+\frac{1}{2}\omega_k,\beta^\vee\rangle\in\nat-\{0\}
}}\Upsilon_{\T^k}(s_\beta(\rho_{\T^k}+\frac{1}{2}\omega_k))=\Upsilon_{\T^k}(s_{e_1+e_{s+1}}(\rho_{\T^k}+\frac{1}{2}\omega_k))+\Upsilon_{\T^k}(s_{e_1+e_{n-s}}(\rho_{\T^k}+\frac{1}{2}\omega_k)).
\end{align*}
However,
$s_{e_1+e_{s+1}}(\rho_{\T^k}+\frac{1}{2}\omega_k)=s_{e_1-e_{s+1}}s_{e_{n-s}-e_n}s_{e_{n-s}+e_n}s_{e_1+e_{n-s}}(\rho_{\T^k}+\frac{1}{2}\omega_k)$
holds.
So, Corollary 2.2.10 implies
\begin{align*}
\sum_{\substack{\beta\in\left(\Delta^{\T^k}\right)^+ \\
 \langle\rho_{\T^k}+\frac{1}{2}\omega_k,\beta^\vee\rangle\in\nat-\{0\}
}}\Upsilon_{\T^k}(s_\beta(\rho_{\T^k}+\frac{1}{2}\omega_k))=0.
\end{align*}
From Jantzen's criterion (Theorem 2.2.11), $M_{\T^k}[\frac{1}{2}]$ is
irreducible.
So, we get (3).

Next, we prove (4).
We put $s=\frac{k-1}{2}$.  
Then, we see
\begin{align*}
\rho_{\T^k}+\frac{1}{2}\omega_k=\sum_{i=1}^{2s+1}\left(s-i+\frac{3}{2}\right)e_i +
 \sum_{j=2s+2}^n \left(n-j+\frac{1}{2}\right)e_j.
\end{align*}
We easily see
\begin{multline*}
\left\{\beta\in \left(\Delta^{\T^k}\right)^+\left|
 \langle\rho_{\T^k}+\frac{1}{2}\omega_k,\beta^\vee\rangle\in\nat-\{0\}\right.\right\}\\
=\{e_i+e_j\mid 1\leqslant i<j\leqslant 2s+1,
 2s+3>i+j\}\\
\cup \left\{e_i\pm e_j\left| 1\leqslant i\leqslant 2s<j\leqslant n,
 \left(s-i+\frac{3}{2}\right)\pm\left(n-j+\frac{1}{2}\right)>0\right. \right\}.
\end{multline*}
Since $3k< 2n+1$, we have $s+\frac{1}{2}\leqslant n-2s-\frac{3}{2}$.
Let $\beta\in \left(\Delta^{\T^k}\right)^+$ be such that 
 $\langle\rho_{\T^k}+\frac{1}{2}\omega_k,\beta^\vee\rangle\in\nat-\{0\}$.
Moreover, we assume $\beta$ is neither $e_1$ nor $e_1+e_{n-s}$.
Then, we can easily see there exists some $\gamma\in\Delta_{\T^k}$ such
that $\langle
s_\beta(\rho_{\T^k}+\frac{1}{2}\omega_k),\gamma^\vee\rangle=0$.
From Corollary 2.2.10, we have
\begin{align*}
\sum_{\substack{\beta\in\left(\Delta^{\T^k}\right)^+ \\
 \langle\rho_{\T^k}+\frac{1}{2}\omega_k,\beta^\vee\rangle\in\nat-\{0\}
}}\Upsilon_{\T^k}(s_\beta(\rho_{\T^k}+\frac{1}{2}\omega_k))=\Upsilon_{\T^k}(s_{e_1}(\rho_{\T^k}+\frac{1}{2}\omega_k))+\Upsilon_{\T^k}(s_{e_1+e_{n-s}}(\rho_{\T^k}+\frac{1}{2}\omega_k)).
\end{align*}
However,
$s_{e_1}(\rho_{\T^k}+\frac{1}{2}\omega_k)=s_{e_{n-s}}s_{e_1+e_{n-s}}(\rho_{\T^k}+\frac{1}{2}\omega_k)$
holds.
So, Corollary 2.2.10 implies
\begin{align*}
\sum_{\substack{\beta\in\left(\Delta^{\T^k}\right)^+ \\
 \langle\rho_{\T^k}+\frac{1}{2}\omega_k,\beta^\vee\rangle\in\nat-\{0\}
}}\Upsilon_{\T^k}(s_\beta(\rho_{\T^k}+\frac{1}{2}\omega_k))=0.
\end{align*}
From Jantzen's criterion (Theorem 2.2.11), $M_{\T^k}[\frac{1}{2}]$ is
irreducible.
So, we get (4).

(5) and (6) is due to [Gyoja 1994] p394. They are proved similarly.
\,\,\,\, Q.E.D.

Lemma 3.6.1 and Lemma 2.2.6 imply:
\begin{lem}
(1) \,\,\,\, Assume that $\gggg=\so(2n+1,\cpx)$ and that $k$ is an even
 positive integer such that $k<n$.
Then, $M_{\T^k}[-t-\frac{1}{2}]\not\subseteq M_{\T^k}[t+\frac{1}{2}]$
 for all $t\in\nat$.

(2) \,\,\,\, Assume that $\gggg=\so(2n+1,\cpx)$ and that $n$ is even.
Then, $M_{\T^k}[-t]\not\subseteq M_{\T^k}[t]$
 for all $t\in\nat$.

(3) \,\,\,\, Assume that $\gggg=\so(2n,\cpx)$ and that $k$ is an even
 positive integer such that $3k<2n$.
Then, $M_{\T^k}[-t-\frac{1}{2}]\not\subseteq M_{\T^k}[t+\frac{1}{2}]$
 for all $t\in\nat$.

(4) \,\,\,\, Assume that $\gggg=\so(2n+1,\cpx)$ and that $k$ is an odd
 positive integer such that $3k<2n+1$.
Then, $M_{\T^k}[-t-\frac{1}{2}]\not\subseteq M_{\T^k}[t+\frac{1}{2}]$
 for all $t\in\nat$.

(5) \,\,\,\, Assume that $\gggg=\sss\ppp(n,\cpx)$ and that $k$ is an even
 positive integer such that $3k>2n$.
Then, $M_{\T^k}[-t-\frac{1}{2}]\not\subseteq M_{\T^k}[t+\frac{1}{2}]$
 for all $t\in\nat$.

(6) \,\,\,\, Assume that $\gggg=\sss\ppp(n,\cpx)$.
Then, $M_{\T^1}[-t]\not\subseteq M_{\T^1}[t]$
 for all $t\in\nat$.
\end{lem}

\subsection{Nonexistence results for the remaining cases}

First, we assume that $\gggg=\so(2n,\cpx)$, $n$ is odd, and $k=n$.
This case is treated in [Boe 1985].
In fact, his result contains:
\begin{lem} (Boe)
Assume $\gggg=\so(2n,\cpx)$ and $n$ is odd.
Then, $M_{\T^n}[-t]\subseteq M_{\T^n}[t]$
 if and only if $t=0$.
\end{lem}

Hereafter we do not consider the above case.
Therefore, from 2.1, we have $w_{\T^k}w_0=w_0w_{\T^k}$. 
The results in 3.3-3.6 and  Corollary 2.2.4 imply that Theorem
3.2.1-3.2.3 is reduced to the following lemma.
\begin{lem} \mbox{}

(1) \,\,\,\,  Assume that $\gggg=\so(2n+1,\cpx)$ and that $k$ is an even
 positive integer such that $3k>2n+1$ and $k<n$.
Then, $M_{\T^k}[-t-1]\not\subseteq M_{\T^k}[t+1]$
 for all $t\in\nat$.

(2) \,\,\,\, Assume that $\gggg=\sss\ppp(n,\cpx)$ and that $k$ is an odd
 positive integer such that $3k\leqslant  2n$ and $1<k$.
Then, $M_{\T^k}[-t-1]\not\subseteq M_{\T^k}[t+1]$
 for all $t\in\nat$. 

(3) \,\,\,\, Assume that $\gggg=\so(2n,\cpx)$ and that $k$ is an odd
 positive integer such that $3k\geqslant 2n$.
Then, $M_{\T^k}[-t-1]\not\subseteq M_{\T^k}[t+1]$
 for all $t\in\nat$.

\end{lem}

\proof
From Lemma 2.2.6, we have only to show that $M_{\T^k}[-1]\not\subseteq
 M_{\T^k}[1]$.

We put $\Omega=W\cdot(\rho_{\T^k}+\omega_k)\cap\PPSK$.
Then, in the settings of (1), (2), and (3) above, $\Omega$ consists of four
elements (say $\lambda_1,...\lambda_4$).
We can write $\lambda_1=\rho_{\T^k}+\omega_k$ and
$\lambda_4=\rho_{\T^k}-\omega_k$.
The remaining two elements are as follows.

If $\gggg=\so(2n+1,\cpx)$ and if $k=2s$ is an even
 positive integer such that $3k>2n+1$ and $k<n$, then we may write:
\begin{gather*}
\lambda_2=
 \left(s+\frac{1}{2}\right)e_1+\sum_{i=2}^{2s}\left(s+\frac{1}{2}-i\right)e_i+\sum_{j=2s+1}^n\left(n+\frac{1}{2}-j\right)e_j, \\
\lambda_3=
 \sum_{i=1}^{2s-1}\left(s+\frac{1}{2}-i\right)e_i+\left(-\frac{1}{2}-s\right)e_{2s}+\sum_{j=2s+1}^n\left(n+\frac{1}{2}-j\right)e_j.
 \end{gather*}

If $\gggg=\sss\ppp(n,\cpx)$ and if  $k=2s+1$ is an odd
 positive integer such that $3k\leqslant  2n$ and $1<k$, then we may write: 
\begin{gather*}
\lambda_2=
 \left(s+1\right)e_1+\sum_{i=2}^{2s+1}\left(s+1-i\right)e_i+\sum_{j=2s+2}^n\left(n+1-j\right)e_j, \\
\lambda_3=
 \sum_{i=1}^{2s}\left(s+1-i\right)e_i+\left(-1-s\right)e_{2s+1}+\sum_{j=2s+2}^n\left(n+1-j\right)e_j.
 \end{gather*}

If $\gggg=\so(2n,\cpx)$ and if $k=2s+1$ is an odd
 positive integer such that $3k\geqslant 2n+1$ and $k<n$, then we may write:
\begin{gather*}
\lambda_2=
 \left(s+1\right)e_1+\sum_{i=2}^{2s+1}\left(s+1-i\right)e_i+\sum_{j=2s+2}^n\left(n-j\right)e_j, \\
\lambda_3=
 \sum_{i=1}^{2s}\left(s+1-i\right)e_i+\left(-1-s\right)e_{2s+1}+\sum_{j=2s+2}^n\left(n-j\right)e_j.
 \end{gather*}

\medskip

{\bf Claim 1} \,\,\,\, {\it $M_{\T^k}(\lambda_3)$ is reducible.}

\medskip

{\it Proof of Claim 1} \,\,\,\,\,\,
We apply Jantzen's Criterion of irreducibility (Theorem 2.2.11).

If $\gggg=\so(2n+1,\cpx)$ and if $k=2s$ is an even
 positive integer such that $3k>2n+1$ and $k<n$, then we have
\begin{align*}
\sum_{\substack{\beta\in\left(\Delta^{\T^k}\right)^+ \\
 \langle\lambda_3,\beta^\vee\rangle\in\nat-\{0\}
}}\Upsilon_{\T^k}(s_\beta\lambda_3) 
 =\Upsilon_{\T^k}(s_{e_1}\lambda_3)
=-[M_{\T^k}(\lambda_4)]\neq 0.
\end{align*}

If $\gggg=\sss\ppp(n,\cpx)$ and if  $k=2s+1$ is an odd
 positive integer such that $3k\leqslant  2n$ and $1<k$, then we may have 
\begin{align*}
\sum_{\substack{\beta\in\left(\Delta^{\T^k}\right)^+ \\
 \langle\lambda_3,\beta^\vee\rangle\in\nat-\{0\}
}}\Upsilon_{\T^k}(s_\beta\lambda_3) &
 =\Upsilon_{\T^k}(s_{e_1}\lambda_3)+
 \Upsilon_{\T^k}(s_{e_1+e_{s+1}}\lambda_3)+ \Upsilon_{\T^k}(s_{e_1+e_{n-s+1}}\lambda_3)\\
&=[M_{\T^k}(\lambda_4)] - [M_{\T^k}(\lambda_4)] - [M_{\T^k}(\lambda_4)]=-[M_{\T^k}(\lambda_4)]\neq 0.
\end{align*}

If $\gggg=\so(2n,\cpx)$ and if $k=2s+1$ is an odd
 positive integer such that $3k\geqslant 2n+1$ and $k<n$, then we have 
\begin{align*}
\sum_{\substack{\beta\in\left(\Delta^{\T^k}\right)^+ \\
 \langle\lambda_3,\beta^\vee\rangle\in\nat-\{0\}
}}\Upsilon_{\T^k}(s_\beta\lambda_3) 
 =\Upsilon_{\T^k}(s_{e_1+e_{n-s+1}}\lambda_3)
=-[M_{\T^k}(\lambda_4)]\neq 0.
\end{align*}
Therefore   we have Claim 1. \,\,\,\, $\Box$

\medskip

{\bf Claim 2} \,\,\,\, There is a non-trivial $U(\gggg)$-homomorphism
$\varphi : M_{\T^k}(\lambda_2)\rightarrow  M_{\T^k}(\lambda_1)$.

\medskip

{\it Proof of Claim 2} \,\,\,\,\,\,
We remark that neither $\lambda_3-\lambda_1$ nor $\lambda_3-\lambda_2$
are contained in $\Q^+$. 
So, Claim 1 implies that $M_{\T^k}(\lambda_4)\subseteq
M_{\T^k}(\lambda_3)$.
On the other hand, Zuckerman's duality theorem ([Boe-Collingwood 1985] 4.9,
[Collingwood-Shelton 1990] Theorem 1.1, [Gyoja 2000]) implies
\begin{align*}
\dim\Hom_{U(\gggg)}(M_{\T^k}(\lambda_4), M_{\T^k}(\lambda_3))=
\dim\Hom_{U(\gggg)}(M_{\T^k}(\lambda_2), M_{\T^k}(\lambda_1)).
\end{align*}
Hence there is a non-trivial $\varphi\in
\Hom_{U(\gggg)}(M_{\T^k}(\lambda_2), (M_{\T^k}(\lambda_1))$.
\,\,\,\, $\Box$

We need:

\medskip

{\bf Claim 3} \,\,\,\, $\Dim(L(\lambda_2))=\dim\nnn_{\T^k}$.

\medskip

Since the multiplicity (the Bernstein degree) of any scalar generalized
Verma module is one, Claim 2 and Claim 3 imply $L(\lambda_2)$ is a
unique irreducible constituent of $M_{\T^k}(\lambda_1)=M_{\T^k}[1]$
whose Gelfand-Kirillov dimension is $\dim \nnn_{\T^k}$.
On the other hand,  $\Dim(M_{\T^k}[-1])=\dim \nnn_{\T^k}$ and
$M_{\T^k}[-1]$ is irreducible.
Hence,  $M_{\T^k}[-1]\not\subseteq M_{\T^k}[1]$ and we have Lemma 3.7.2.
So, we have only to show Claim 3.

{\it Proof of Claim 3} \,\,\,\,\,\,

Let $\sigma\in W$ be the longest element (with respect to the length
$\ell(\cdot)$) in 
\begin{align*}
\{w\in W\mid \mbox{$\lambda_2$ is dominant with
respect to $w\Delta^+$}\}.
\end{align*}

Then, $\sigma\rho$ and $\lambda_2$ are contained in the same closed Weyl
chamber and we can regard $L(\lambda_2)$ as a limit of $L(\sigma\rho)$.
Namely, $T^{\lambda_2}_{\sigma\rho}(L(\sigma\rho))\cong L(\lambda_2)$.
Since $\Dim(L(\sigma\rho))=\Dim(\lambda_2)$, Claim 3 is reduced to the
following Claim 4. (Cf.\ Lemma 2.1.6.)

\medskip

{\bf Claim 4} \,\,\,\, $\sigma\lreq w_{\T^k}w_0$.

\medskip

{\it Proof of Claim 4} \,\,\,\,\,\,
Claim 4 is obtained by the algorithm described in [Barbasch-Vogan 1982]
p171-175.

First, we assume  $\gggg=\sss\ppp(n,\cpx)$ and that $k$ is an odd
 positive integer such that $3k\leqslant  2n$ and $1<k$.

We  put $\delta=10^{-23}$. 
(In fact, $\delta$ can be any real number such that $0<\delta<\frac{1}{2}$.)
Put $\lambda_2^\prime= \lambda_2-\delta\omega_k$.
We also put $a_i=\langle e_i,\lambda_2^\prime\rangle$ for $1\leqslant
i\leqslant n$. 
Namely, $\lambda^\prime_2=\sum_{i=1}^n a_ie_i$.
Then, we easily see $|a_1|,....,|a_n|$ are distinct.
Let $b_1,...,b_n$ be positive real numbers such that
$b_1>b_2>\cdots>b_n$ and $\{b_1,..., b_n\}=\{|a_1|,....,|a_n|\}$.
So, there is a permutation $\tau\in\Sth_n$ such that $|a_i|=b_{\tau(i)}$
for $1\leqslant i\leqslant n$. For $1\leqslant i\leqslant n$, we put
$c_i=\frac{a_i}{|a_i|}(n-\tau(i)+1)$.
Then, we have $\sigma\rho=\sum_{i=1}^n c_n e_n$.
Barbasch and Vogan attach $\sigma$ to the sequence
$(c_1,...,c_n,-c_n,...,-c_1)$ ([Barbasch-Vogan 1982] p173).
Applying the Robinson-Schensted algorithm to this sequence, we get a pair
of Young tableaux. (These Young tableaux have the same shape.)
Remark that the Robinson-Schemsted algorithm in [Barbasch-Vogan 1982] is
a little bit
different from the usual one (see [Barbasch-Vogan 1982] p171).
For our purpose, the important information is the shape of these Young
tableaux.
In order to obtain such a Young diagram, we need not compute
$c_1,...,c_n$.
In fact, applying  the Robinson-Schensted algorithm to the sequence
$(a_1,....,a_n,-a_n,...,-a_1)$ directly, we have the same Young diagram.
In this case, the Young diagram corresponds to the partition of $2n$,
$(2n-2k, k+1, k-1)$.
It corresponds to the symbol
$\left(\begin{array}{rrr}\frac{k-1}{2} & & n-k+1 \\  & \frac{k+1}{2} &\end{array}
 \right)$.

On the other hand, $w_{\T^k}w_0$ corresponds to the partition of $2n$,
$(2n-2k, k,  k)$.
It corresponds to the symbol
$\left(\begin{array}{rrr}\frac{k+1}{2} & & n-k+1 \\  & \frac{k-1}{2} & \end{array}
\right)$.
Hence, from [Barbasch-Vogan 1982] Theorem 18, we have $\sigma\lreq w_{\T^k}w_0$.

Next, we consider the case of $\gggg=\so(2n+1,\cpx)$.
We define $\lambda_2^\prime$, $a_i$ \,\, $(1\leqslant i\leqslant n)$,
$b_i$ \,\, $(1\leqslant i\leqslant n)$, $\tau\in\Sth_n$, and $c_i$ \,\,
$(1\leqslant i\leqslant n)$ in the same way as the case of
$\gggg=\sss\ppp(n,\cpx)$.
In the Barbasch-Vogan algorithm, $W$ is regarded as the Weyl group of
the type $C_n$ rather than $B_n$.
So, $\sigma$ is attached to $(c_1,...,c_n,-c_n,...,-c_1)$.
Again, we may apply the Robinson-Schensted algorithm to the sequence
$(a_1,...,a_n,-a_n,...,-a_1)$ directly and obtain a partition $(k+1, k-1,
2n-2k)$ of $2n$.
It corresponds to a symbol
$\left(\begin{array}{rrr}n-k & & \frac{k}{2}\\  & \frac{k}{2}+1 &\end{array}
 \right)$.
On the other hand, $w_{\T^k}w_0$ corresponds to the partition of $2n$,
$(k,  k, 2n-2k)$.
It corresponds to a symbol
$\left(\begin{array}{rrr} n-k & & \frac{k}{2}+1 \\  & \frac{k}{2} & \end{array}
\right)$.
Hence, from [Barbasch-Vogan 1982] Theorem 18, we have $\sigma\lreq w_{\T^k}w_0$.

Finally, we consider the case of $\gggg=\so(2n,\cpx)$.
In this case, we put $\lambda_2^\prime=
\lambda_2-\delta\omega_k+\frac{\delta}{2}e_n$.
Here, $\delta$ is a fixed real number such that $0<\delta<\frac{1}{2}$.
We define $a_i$ \,\, $(1\leqslant i\leqslant n)$,
$b_i$ \,\, $(1\leqslant i\leqslant n)$, $\tau\in\Sth_n$, and $c_i$ \,\,
$(1\leqslant i\leqslant n)$ in the same way as the case of
$\gggg=\sss\ppp(n,\cpx)$.
$\sigma$ is attached to $(c_1,...,c_n,-c_n,...,-c_1)$.
Again, we may apply the Robinson-Schensted algorithm to the sequence
$(a_1,...,a_n,-a_n,...,-a_1)$ directly and obtain a partition $(k+1, k-1,
2n-2k)$ of $2n$.
It corresponds to a symbol
$\left(\begin{array}{cc} n-k &  \frac{k+3}{2} \\  0  & \frac{k+1}{2}  \end{array}
\right)$.
On the other hand, $w_{\T^k}w_0$ corresponds to the partition of $2n$,
$(k,  k, 2n-2k)$.
It corresponds to a symbol
$\left(\begin{array}{cc} n-k &  \frac{k+1}{2} \\  0  & \frac{k+3}{2} \end{array}
\right)$.
Hence, from [Barbasch-Vogan 1982] Theorem 18, we have $\sigma\lreq w_{\T^k}w_0$.
\,\,\,\, Q.E.D.


\setcounter{section}{4}
\setcounter{subsection}{0}

\section*{\S\,\, 4.\,\,\,\, Exceptional algebras}

As in \S 3, we write $\T^k=\Pi-\{\alpha_k\}$,
$\omega_k=\omega_{\alpha_k}$, $c_k=c_{\alpha_k}$, $d_k=d_{\alpha_k}$, etc.

\subsection{Even parabolic subalgebras}

Let $u\in\gggg$ be a nilpotent element.
From the Jacobson-Morozov theorem, there is an $\sss\llll_2$-triple $(v,
h, u)$ in $\gggg$.
Namely, we have $[u,v]=h$, $[h,u]=2u$, and $[h,v]=-2v$.
It is well known that any eigenvalue of the operator
$\ad(h)\in\End(\gggg)$ is an integer. 
$u$ is called an even nilpotent element  if each eigenvalue of $\ad(h)$
is an even number.
Put $\ppp_u=\sum_{i\geqslant 0}\{X\in\gggg\mid [h, X] =iX\}$.
\begin{dfn}
A parabolic subalgebra $\ppp$ is called even, if there exists an even
 nilpotent element $u$ such that $\ppp=\ppp_u$.
\end{dfn}

For an even parabolic subalgebra $\ppp_u$, the Richardson orbit
$\hol_{\ppp_u}$ contains $u$.

The following result is well-known. For example, it is an easy consequence
of  [Hesselink 1978] p218
and [Yamashita 1986] Lemma 3.5.
\begin{lem}
Let $\ppp_{\T}$ be an even parabolic subalgebra of $\gggg$.
Then,  the moment map $m_{\T} : T^\ast X_{\T}\rightarrow \overline{\hol_{\ppp_{\T}}}$ is birational.
\end{lem}

\subsection{${\bf G_2}$}

In the regular integral case, homomorphisms between (not necessarily scalar)
generalized Verma modules are classified for $G_2$ in  [Boe-Collingwood
1990].

I imagine in the case of $G_2$ the classification of the homomorphisms
between scalar generalized Verma modules is well-known, but I would like to state the result
for the completeness.

Let $\gggg$ be a simple Lie algebra of the type $G_2$.
Then, we may write $\Pi=\{\alpha_1,\alpha_2\}$ such that
$\langle\alpha_1,\alpha_1\rangle>\langle\alpha_2,\alpha_2\rangle $.

\begin{thm}
Let $\gggg$ be a simple Lie algebra of the type $G_2$ and let $k$ be $1$ or $2$.
Then, $M_{\T^k}[-t]\subseteq M_{\T^k}[t]$ if and only if $t\in\nat$.
\end{thm}

\proof
We see $c_1=c_2=\frac{1}{2}$, $d_1=\frac{3}{2}$, and
$d_2=\frac{5}{2}$.
Thus, for $k=1,2$, $\rho_{\T^k}+t\omega_k$ is integral if and only if
$t-\frac{1}{2}\in\itg$.
For $k=1,2$, it is easy to see that $w_{\T^k}w_0=w_0w_{\T^k}$ but they
are not a Duflo involution of $W$.
Moreover, we have:
\begin{lem} 
Let $k$ be $1$ or $2$.
If $\rho_{\T^k}+t\omega_k$ is singular and integral, then $M_{\T^k}[t]$
 is irreducible.
In particular, $M_{\T^k}[\frac{1}{2}]$ is irreducible.
\end{lem}

\proof 
$\gggg$ has only three special nilpotent orbits, namely the regular
nilpotent orbit, the subregular nilpotent orbit, and $\{0\}$.
(Cf.\ [Carter 1985])
Hence, the Gelfand-Kirillov dimension of  any infinite-dimensional
irreducible constituent of $M_{\T^k}[t]$ is
$\dim\nnn_{\T^k}=\Dim(M_{\T^k}[t])$, if $\rho_{\T^k}+t\omega_k$ is integral.
On the other hand, the multiplicity (the Bernstein degree) of
$M_{\T^k}[t]$ is one.
So, if  $\rho_{\T^k}+t\omega_k$ is integral and singular, $M_{\T^k}[t]$
is irreducible. \,\,\,\, $\Box$

We continue the proof of Theorem 4.2.1.
From Lemma 4.2.2 and  Lemma 2.2.6, we have
$M_{\T^k}[-\frac{1}{2}-t]\not\subseteq M_{\T^k}[\frac{1}{2}+t]$ for all
$t\in \nat$.

For $t\in\itg$, $\rho_{\T^k}+t\omega_k$ is half-integral and its
integral Weyl group is of type $A_2\times A_2$.
This case, we easily see $w_{\T^k}w_0$ is a Duflo involution of the
integral Weyl group.
Since $\rho_{\T^k}+ t\omega_k$ is dominant and  regular if $t\in\nat-\{0\}$, 
Theorem 1.4.2 implies $M_{\T^k}[-1-t]\subseteq M_{\T^k}[1+t]$ for all
$t\in\nat$.
\,\,\,\, Q.E.D.

\subsection{${\bf F_4}$}

For the simple algebra of the type $F_4$, [Boe-Collingwood] treated
 the regular
 integral case.
The half-integral case is somewhat easier, but I would like to mention
 the results for the completeness.

We consider the root system $\Delta$ for a simple Lie
algebra $\gggg$ of the type $F_4$. (For example, see [Knapp 2002] p691.)
We can choose an orthonormal basis $e_1,...,e_4$ of $\hhh^\ast$ such that 
\begin{align*}
\Delta = \{ \pm e_i\pm e_j\mid 1\leqslant i< j\leqslant
4\}\cup
\{\pm e_i\mid 1\leqslant i\leqslant 4\}
\cup\left\{\frac{1}{2}(\pm e_1\pm e_2\pm e_3 \pm e_4 )\right\}.
\end{align*}
We choose a positive system as follows.
\begin{align*}
\Delta^+ = \{ e_i\pm e_j\mid 1\leqslant i< j\leqslant
4\}\cup
\{e_i\mid 1\leqslant i\leqslant 4\}\cup
\left\{\frac{1}{2}(e_1\pm e_2\pm e_3 \pm e_4 )\right\}.
\end{align*}
Put $\alpha_1=\frac{1}{2}(e_1-e_2-e_3-e_4 )$, $\alpha_2=e_4$,
$\alpha_3=e_3-e_4$, and $\alpha_4=e_2-e_3$.
Then, $\Pi=\{\alpha_1,...,\alpha_4\}$.
\[
 \begin{array}{ccccccc} 1 & - & 2 & \Leftarrow & 3 & - & 4
\end{array}
\]
Since $\PP=\Q$ in this case, we have $c_k=\frac{1}{2}$ for $1\leqslant
k\leqslant 4$.

The result is:

\begin{thm} ((1) is due to [Lepowsky 1975a].)\mbox{} 

(1) \,\,\,\,  
$M_{\T^1}[-t]\subseteq M_{\T^1}[t]$ if and only if $t\in\frac{1}{2}\nat$.

(2) \,\,\,\, 
$M_{\T^2}[-t]\subseteq M_{\T^2}[t]$ if and only if $t\in\nat$.

(3) \,\,\,\, 
$M_{\T^3}[-t]\subseteq M_{\T^3}[t]$ if and only if $t\in\nat$.

(4) \,\,\,\, 
$M_{\T^3}[-t]\subseteq M_{\T^3}[t]$ if and only if $t\in\frac{1}{2}\nat$.

\end{thm}

\proof

(1) is proved in [Lepowsky 1975a] (Theorem 1.2).
So, we consider the other cases.

Since $c_k=\frac{1}{2}$ for $k=1,2,3,4$, $M_{\T^k}[-t]\subseteq M_{\T^k}[t]$
implies $t\in\frac{1}{2}\nat$.

For each  $k=1,3,4$, $\ppp_{\T^k}$ is  an even parabolic subalgebra.
(Cf.\ [Carter 1985] p401.)
Hence, $M_{\T^k}[-t]\subseteq M_{\T^k}[t]$ for $t\in \nat$ for
$k=1,3,4$.

We consider the case of  $k=2$.
$\rho_{\T^2}+t\omega_2$ is half-integral for $t\in\itg$.
The integral root system $\Delta_{\rho_{\T^2}}$ for $\rho_{\T^2}+t\omega_2$ \,\,\, $(t\in\nat)$ is of
the type $A_1\times B_3$.
$\T^2$ corresponds to the set of the long simple roots in the $B_2$
factor of $\Delta_{\rho_{\T^2}}$ and the simple root in the $A_1$-factor.
Hence, Theorem 3.2.1 (2a), Lemma 1.4.4, and Theorem 1.4.2 imply that
$w_{\T^2}w_0$ is a Duflo involution in the integral Weyl group.
So, Theorem 1.4.2 implies $M_{\T^2}[-t]\subseteq M_{\T^2}[t]$ for $t\in
\nat$.

For $k=2,3$, we can prove that $M_{\T^k}[\frac{1}{2}]$ is irreducible
using Jantzen's criterion (Theorem 2.2.11).
Hence, Lemma 2.2.6 implies that $M_{\T^k}[-\frac{1}{2}-t]\not\subseteq M_{\T^k}[\frac{1}{2}+t]$
for $t\in\nat$ and $k=2,3$.

For $k=4$, $\rho_{\T^k}+(\frac{1}{2}+t)\omega_4$ is half-integral for
$t\in\nat$.
The integral root system $\Delta_{\rho_{\T^2}}$ for $\rho_{\T^2}+t\omega_2$ \,\,\, $(t\in\nat)$ is of
the type $A_1\times B_3$.
In this case $w_{\T^4}$ is the non-trivial element of the $A_1$-factor
of the integral Weyl group.
So, it is a Duflo involution in the integral Weyl group and we have
$M_{\T^4}[-\frac{1}{2}-t]\subseteq M_{\T^4}[\frac{1}{2}+t]$
for $t\in\nat$.
\,\,\,\, Q.E.D.  

\subsection{${\bf E_6}$}

We consider the root system $\Delta$ for a simple Lie
algebra $\gggg$ of the type $E_6$. 
Put $\kappa=\frac{1}{2\sqrt{3}}$.
We can choose an orthonormal basis $e_1,...,e_6$ of $\hhh^\ast$ such that 
\begin{multline*}
\Delta = \{ e_i -e_j\mid 1\leqslant i, j\leqslant
6, i\neq j\} \\
\cup 
\left\{\pm\sum_{i=1}^6\left(\kappa+(-1)^{n(i)}\frac{1}{2}\right)e_i\left|
 \mbox{$n(i)=1$ or $n(i)=0$ for $1\leqslant i\leqslant 6$},
 \sum_{i=1}^6n(i)=3\right.  \right\}\\
\cup\left\{\pm 2\kappa\sum_{i=1}^6e_i \right\}.
\end{multline*}

We choose a positive system as follows.
\begin{multline*}
\Delta^+ = \{ e_i -e_j\mid 1\leqslant i< j\leqslant
6 \} \\
\cup 
\left\{\sum_{i=1}^6\left(\kappa+(-1)^{n(i)}\frac{1}{2}\right)e_i\left|
 \mbox{$n(i)=1$ or $n(i)=0$ for $1\leqslant i\leqslant 6$},
 \sum_{i=1}^6n(i)=3\right.  \right\}\\
\cup\left\{2\kappa\sum_{i=1}^6e_i \right\}.
\end{multline*}
Put $\alpha_i=e_i-e_{i+1}$  \,\,\, $(1\leqslant i\leqslant 5)$ and 
$\alpha_6=\sum_{i=1}^3\left(\kappa-\frac{1}{2}\right)e_i+\sum_{i=4}^6\left(\kappa+\frac{1}{2}\right)e_i$.
Then, $\Pi=\{\alpha_1,...,\alpha_6\}$.
\[
 \begin{array}{ccccccccc} 1 & - & 2 & - & 3 & - & 4 & - & 5 \\
&  & & & \mid & & & & \\
&  & & & 6 & & & & 
\end{array}
\]

In this case, $w_0$ is not contained in the center of the Weyl group.
So, $w_{\T}w_0=w_0w_{\T}$ may fail for some  $\T\subseteq\Pi$.
In fact, $w_{\T^k}w_0=w_0w_{\T^k}$ holds for $k=3,6$, but it fails for $k=1,2,4,5$.

Since the Dynkin diagram of $E_6$ has a symmetry, the cases of $\T^4$
and $\T^5$ are similar to $\T^2$ and $\T^1$, respectively.
So, we only consider $\T^1, \T^2, \T^3,$ and $\T^6$.

We have:
\begin{thm} ((1) is due to [Boe 1985].) \mbox{}
 
(1) \,\,\,\,  
$M_{\T^1}[-t]\subseteq M_{\T^1}[t]$ if and only if $t=0$.

(2) \,\,\,\, 
$M_{\T^2}[-t]\subseteq M_{\T^2}[t]$ if and only if $t=0$.

(3) \,\,\,\, 
$M_{\T^3}[-t]\subseteq M_{\T^3}[t]$ if and only if $t\in\nat$.

(4) \,\,\,\, 
$M_{\T^6}[-t]\subseteq M_{\T^3}[t]$ if and only if $t\in\nat$.

\end{thm}

\proof

(1) is proved in [Boe 1985].
So, we consider the other cases.

First, we prove (2).
In this case, we have $\rho^{\T^2}=3\omega_2$.
We also have $\rho^{\T^2}\in\Q$ and $c_2=\frac{3}{2}$.
Hence, $M_{\T^2}[-t]\subseteq M_{\T^2}[t]$ implies $t\in\frac{3}{2}\nat
$.
If $t\in 3\nat$, then $\rho_{\T^2}+t\omega_2$ is integral.

Assume that $t-\frac{3}{2}\in 3\nat$. 
Then, the integral root system of
$\rho_{\T^2}+t\omega_2$ is of the type $A_1\times A_5$.
The set of simple roots is
$\{\alpha_1,\alpha_6, \alpha_3,\alpha_4,\alpha_5,\beta\}$.
Here,
$\beta=\alpha_1+2\alpha_2+2\alpha_3+\alpha_4+\alpha_6=\sum_{i=1}^2\left(\kappa+\frac{1}{2}\right)e_i+\sum_{i=3}^5\left(\kappa-\frac{1}{2}\right)e_i+\left(\kappa+\frac{1}{2}\right)e_6$.
Among them, $\alpha_1$ is the simple root in the $A_1$-factor.

For $t\in\frac{3}{2}\nat$, we denote by $W_{(t)}$ the integral Weyl group
of $\rho_{\T^2}+t\omega_3$.
Let $w^t_0$ be the longest element of the integral Weyl group.
We easily see that $w_{\T^2}w^t_0$ is not an involution.
We also see $\rho_{\T^2}+t\omega_3$ is  dominant regular for
$t\in\frac{3}{2}\nat-\{0\}$.
Since $w_{\T^2}(\rho_{\T^2}-t\omega_3)=-(\rho_{\T^2}+t\omega_3)$,
one of the following two conditions must hold.

(a) \,\,\,\, There is no $w\in W_{(t)}$ such that $w(\rho_{\T^2}+t\omega_3)=\rho_{\T^2}-t\omega_3$.

(b) \,\,\,\,
$w_{\T^2}w^t_0(\rho_{\T^2}+t\omega_3)=\rho_{\T^2}-t\omega_3$.

If (a) holds, then $M_{\T^2}[t]$ has an infinitesimal character
different from that of  $M_{\T^2}[-t]$.
Hence, we have $M_{\T^2}[-t]\not\subseteq M_{\T^2}[t]$.

Assume (b) holds.  Then, Theorem 1.4.2 implies that
$M_{\T^2}[-t]\not\subseteq M_{\T^2}[t]$, since $w_{\T^2}w^t_0$ is not an
involution.

Next, we prove (3) and (4).
Let $k$ be either $3$ or $6$.
Then we have $w_{\T^k}w_0=w_0w_{\T^k}$ and $\ppp_k$ is an even
parabolic subalgebra (see [Carter 1985] p402).
Hence, Lemma 4.1.2 and Corollary 2.2.9 imply $M_{\T^k}[-t]\subseteq
M_{\T^k}[t]$ for all $t\in\nat$.
Since $c_k=\frac{1}{2}$ in these case, we have only to show
$M_{\T^2}[-\frac{1}{2}-t]\not\subseteq M_{\T^2}[\frac{1}{2}+t]$ for all
$t\in\nat$.
In this case,  we can check the irreducibility of $M_{\T^k}[\frac{1}{2}]$
 via Jantzen's criterion (Theorem 2.2.11).
So, from Lemma 2.2.6, we have the desired result.
We describe the computation  briefly.

First, we consider the case of $k=3$.
In this case $d_3=\frac{7}{2}$ and $\rho_{\T^3}+\frac{1}{2}\omega_3$ is integral.
We put $\gamma_1=\sum_{i=1}^62\kappa e_i$ and  $\gamma_2=\sum_{i=1}^3\left(\kappa+\frac{1}{2}\right)e_i+\sum_{i=4}^6\left(\kappa-\frac{1}{2}\right)e_i$.

Put $\Xi=\left\{\beta\in\left. \left(\Delta^{\T^3}\right)^+ \right|
 \langle\rho_{\T^3}+\frac{1}{2}\omega_3,\beta^\vee\rangle\in\nat-\{0\}\right\}$.
Then, $\gamma_1,\gamma_{2}\in\Xi$.
For $\beta\in\Xi-\{\gamma_1, \gamma_2\}$, we can find
 $\eta\in\Delta_{\T^3}$ such that $\langle\beta,\eta\rangle=0$.
Hence, we have
 $\Upsilon_{\T^3}(s_\beta(\rho_{\T^3}+\frac{1}{2}\omega_3))=0$.
Moreover, we have
$\Upsilon_{\T^3}(s_{\gamma_{1}}(\rho_{\T^3}+\frac{1}{2}\omega_3))
=-\Upsilon_{\T^3}(\rho_{\T^3}-\frac{1}{2}\omega_3)$ and
$\Upsilon_{\T^3}(s_{\gamma_{2}}(\rho_{\T^3}+\frac{1}{2}\omega_3)) =\Upsilon_{\T^3}(\rho_{\T^3}-\frac{1}{2}\omega_3)$.
Hence, we have
\begin{align*}
\sum_{\beta\in\Xi}\Upsilon_{\T^3}(s_{\beta}(\rho_{\T^3}+\frac{1}{2}\omega_3))=0.
\end{align*}
This means that Jantzen's criterion is satisfied.

Finally,  we consider the case of $k=6$.
In this case $d_6=\frac{11}{2}$ and $\rho_{\T^6}+\frac{1}{2}\omega_6$ is integral.
We put $\gamma_1=\sum_{i=1}^62\kappa e_i$ and  $\gamma_2=\sum_{i=1}^6\left(\kappa+(-1)^{i-1}\frac{1}{2}\right)e_i$.

Put $\Xi=\left\{\beta\in\left. \left(\Delta^{\T^6}\right)^+ \right|
 \langle\rho_{\T^6}+\frac{1}{2}\omega_6,\beta^\vee\rangle\in\nat-\{0\}\right\}$.
Then, $\gamma_1,\gamma_{2}\in\Xi$.
For $\beta\in\Xi-\{\gamma_1, \gamma_2\}$, we can find
 $\eta\in\Delta_{\T^6}$ such that $\langle\beta,\eta\rangle=0$.
Hence, we have
 $\Upsilon_{\T^6}(s_\beta(\rho_{\T^6}+\frac{1}{2}\omega_6))=0$.
Moreover, we have
$\Upsilon_{\T^6}(s_{\gamma_{1}}(\rho_{\T^6}+\frac{1}{2}\omega_6))
=\Upsilon_{\T^6}(\rho_{\T^3}-\frac{1}{2}\omega_6)$ and
$\Upsilon_{\T^6}(s_{\gamma_{2}}(\rho_{\T^6}+\frac{1}{2}\omega_6)) =-\Upsilon_{\T^6}(\rho_{\T^6}-\frac{1}{2}\omega_6)$.
Hence, we have
\begin{align*}
\sum_{\beta\in\Xi}\Upsilon_{\T^6}(s_{\beta}(\rho_{\T^6}+\frac{1}{2}\omega_6))=0.
\end{align*}
This means that Jantzen's criterion is satisfied.
\,\,\,\, Q.E.D.

{\bf Remark} \,\,\,\, In the case of $k=6$, the non-existence of the homomorphism is proved in
[Boe-Collingwood 1990] for the regular integral case.

\subsection{${\bf E_7}$}

Let $\gggg$ be a simple Lie algebra of the type $E_7$.

We fix an orthonormal basis $e_1,...,e_8$ in $\rel^8$.
We identify $\hhh^\ast$ with $\{v\in\rel^8\mid \langle v,
e_1-e_2\rangle=0\}$ so that
\begin{multline*}
\Delta=\{\pm (e_1+e_2)\}\cup\{\pm e_i\pm e_j\mid 3\leqslant i<j\leqslant
8\}\\
\cup \left\{\pm\frac{1}{2}\left(e_1+e_2+\sum_{i=3}^8(-1)^{n(i)}e_i\right)\left|
\mbox{$n(i)$ is either $0$ or $1$ for $3\leqslant i\leqslant 8$ and
 $\sum_{i=3}^8n(i)$ is even.} \right.\right\}
\end{multline*}

We choose a positive system as follows.
\begin{multline*}
\Delta^+=\{ (e_1+e_2)\}\cup\{ e_i\pm e_j\mid 3\leqslant i<j\leqslant
8\}\\
\cup \left\{\frac{1}{2}\left(e_1+e_2+\sum_{i=3}^8(-1)^{n(i)}e_i\right)\left|
\mbox{$n(i)$ is either $0$ or $1$ for $3\leqslant i\leqslant 8$ and
 $\sum_{i=3}^8n(i)$ is even.} \right.\right\}
\end{multline*}

Put $\alpha_i=e_{i+2}-e_{i+3}$ for $1\leqslant i\leqslant 5$,
$\alpha_6=e_7+e_8$, and
$\alpha_7=\frac{1}{2}(e_1+e_2-e_3-e_4-e_5-e_6-e_7-e_8)$.
Then, $\Pi=\{\alpha_1,...,\alpha_7\}$ is the set of simple roots in
$\Delta^+$.

\[
 \begin{array}{ccccccccccc} 1 & - & 2 & - & 3 & - & 4 & - & 6 & - & 7\\
& &  & & & & \mid & & & & \\
& &  & & & & 5 & & & & 
\end{array}
\]
We have:
\begin{thm} (The case of $k=1$ is due to [Boe 1985].) \mbox{}

(1) \,\,\,\,  Assume $k\in\{1,3,5,6,7\}$.
Then, $M_{\T^k}[-t]\subseteq M_{\T^k}[t]$ if and only if $t\in\nat$.

(2) \,\,\,\,  Assume $k\in\{2,4\}$.
Then, $M_{\T^k}[-t]\subseteq M_{\T^k}[t]$ if and only if $t\in\frac{1}{2}\nat$.

\end{thm}

\proof

For the simple Lie algebra of the type $E_7$, $w_0$ is contained in the
center of the Weyl group.
Moreover, all the maximal parabolic
subalgebras are even. (See [Carter 1985] p403-404.)
So, Lemma 4.1.2 and Corollary 2.2.9 imply that, for all $1\leqslant k\leqslant 8$, $t\in\nat$ implies $M_{\T^k}[-t]\subseteq M_{\T^k}[t]$.
Hence, we have only to take care of the case of $t-\frac{1}{2}\in\nat$.

The case of $k=1$ is due to [Boe 1985].
In fact  $M_{\T^1}[-t]\not\subseteq M_{\T^1}[t]$ for $t-\frac{1}{2}\in\nat$.

For the case $k=2$, we have $d_2=7$ and $\omega_2\in\Q$.
The integral Weyl group for $\rho_{\T^2}+t\omega_2$ \,\,\,
$(t\in\frac{1}{2}+\nat)$
is of the type $A_1\times D_6$.
In fact, we have $\Pi_{\rho_{\T^2}+t\omega_2}=\{ e_3+e_4\}\cup \T^2$.
From Theorem 1.4.2 and Theorem 3.2.3, we have  $w_{\T^2}w_0$ is a Duflo
involution of the integral Weyl group for $\rho_{\T^2}+t\omega_2$.
So, Theorem  1.4.2, Lemma 2.2.6 imply that $M_{\T^2}[-t]\subseteq
M_{\T^2}[t]$ for $t-\frac{1}{2}\in\nat$.

For the case $k=3$, we have $d_3=5$ and $\omega_3=\frac{3}{2}(e_1+e_2)+e_3+e_4+e_5\not\in\Q$.
So, we have $c_3=1$.
Hence, $M_{\T^3}[-t]\not\subseteq
M_{\T^3}[t]$ for $t-\frac{1}{2}\in\nat$.

For the case $k=4$, we have $d_4=4$ and $\omega_4\in\Q$.
The integral Weyl group for $\rho_{\T^4}+t\omega_4$ \,\,\,
$(t\in\frac{1}{2}+\nat)$
is of the type $A_1\times D_6$.
In fact, we have $\Pi_{\rho_{\T^4}+t\omega_4}=\{ e_5+e_6\}\cup \T^4$.
From Theorem 1.4.2 and Theorem 3.2.3, we have  $w_{\T^4}w_0$ is a Duflo
involution of the integral Weyl group for $\rho_{\T^4}+t\omega_4$.
So, Theorem  1.4.2, Lemma 2.2.6 imply that $M_{\T^4}[-t]\subseteq
M_{\T^4}[t]$ for $t-\frac{1}{2}\in\nat$.

For the case $k=5$, we have $d_5=7$ and $\omega_7=e_1+e_2+\frac{1}{2}\left(e_3+e_4+e_5+e_6+e_7-e_8\right)\not\in\Q$.
So, we have $c_5=1$.
Hence, $M_{\T^5}[-t]\not\subseteq
M_{\T^5}[t]$ for $t-\frac{1}{2}\in\nat$.

For the cases $k=6$ and $k=7$, we have $d_6=\frac{11}{2}$ and $d_7=\frac{17}{2}$.
In this case, we can show the irreducibility of $M_{\T^k}[\frac{1}{2}]$
\,\, $(k=6,7)$
via Jantzen's criterion (Theorem 2.2.11).
So, from Lemma 2.2.6, for $k=6,7$, we have $M_{\T^k}[-t]\not\subseteq
M_{\T^k}[t]$ for $t-\frac{1}{2}\in\nat$.
We describe the computation briefly.

First, we consider the case of $k=6$.
We remark that $\PP_{\T^k}^{++}\cap W(\rho_{\T^6}+\frac{1}{2}\omega_6)$
consists of the $3$ elements (say $\lambda_1,\lambda_2, \lambda_3$).
We put  $\lambda_1=\rho_{\T^6}+\frac{1}{2}\omega_6$ and $\lambda_3=\rho_{\T^6}-\frac{1}{2}\omega_6$.
The remaining element is $\lambda_2=\frac{1}{2}e_1+\frac{1}{2}e_2+2e_3-e_4-e_6-2e_7-3e_8$.
We put $\gamma_1=e_1+e_2$, $\gamma_2=-\frac{1}{2}e_1-\frac{1}{2}e_2+e_3+2e_5-2e_6-3e_7-e_8$, $\gamma_3=-\frac{1}{2}e_1-\frac{1}{2}e_2+e_3+3e_4+2e_5-2e_6-e_8$, and $\gamma_4=\frac{1}{2}e_1+\frac{1}{2}e_2+2e_3+e_4-e_6-2e_7-3e_8$.
Put $\Xi=\left\{\beta\in\left. \left(\Delta^{\T^6}\right)^+ \right|
 \langle\rho_{\T^6}+\frac{1}{2}\omega_6,\beta^\vee\rangle\in\nat-\{0\}\right\}$.
Then, $\gamma_1,...,\gamma_{4}\in\Xi$.
For $\beta\in\Xi-\{\gamma_1,...,\gamma_{4}\}$, we can find
 $\eta\in\Delta_{\T^6}$ such that $\langle\beta,\eta\rangle=0$.
Hence, we have
 $\Upsilon_{\T^6}(s_\beta(\rho_{\T^6}+\frac{1}{2}\omega_6))=0$.
Moreover, we have
$\Upsilon_{\T^6}(s_{\gamma_{1}}(\rho_{\T^6}+\frac{1}{2}\omega_6)) =-\Upsilon_{\T^6}(\lambda_3)$,
$\Upsilon_{\T^6}(s_{\gamma_{2}}(\rho_{\T^6}+\frac{1}{2}\omega_6)) =\Upsilon_{\T^6}(\lambda_3)$,
$\Upsilon_{\T^6}(s_{\gamma_{3}}(\rho_{\T^6}+\frac{1}{2}\omega_6))
=-\Upsilon_{\T^6}(\lambda_2)$, and
$\Upsilon_{\T^6}(s_{\gamma_{4}}(\rho_{\T^6}+\frac{1}{2}\omega_6))
 =\Upsilon_{\T^6}(\lambda_2)$.
Hence, we have
\begin{align*}
\sum_{\beta\in\Xi}\Upsilon_{\T^6}(s_{\beta}(\rho_{\T^6}+\frac{1}{2}\omega_6))=0.
\end{align*}
This means that Jantzen's criterion is satisfied.

First, we consider the case of $k=7$.
We put $\gamma_1=e_1+e_2$ and  $\gamma_2=\frac{1}{2}(e_1+e_2+e_3-e_4+e_5-e_6+e_7+e_8)$.

Put $\Xi=\left\{\beta\in\left. \left(\Delta^{\T^8}\right)^+ \right|
 \langle\rho_{\T^8}+\frac{1}{2}\omega_8,\beta^\vee\rangle\in\nat-\{0\}\right\}$.
Then, $\gamma_1,\gamma_{2}\in\Xi$.
For $\beta\in\Xi-\{\gamma_1, \gamma_2\}$, we can find
 $\eta\in\Delta_{\T^8}$ such that $\langle\beta,\eta\rangle=0$.
Hence, we have
 $\Upsilon_{\T^8}(s_\beta(\rho_{\T^8}+\frac{1}{2}\omega_8))=0$.
Moreover, we have
$\Upsilon_{\T^8}(s_{\gamma_{1}}(\rho_{\T^8}+\frac{1}{2}\omega_8))
=\Upsilon_{\T^8}(\rho_{\T^8}-\frac{1}{2}\omega_8)$ and
$\Upsilon_{\T^8}(s_{\gamma_{2}}(\rho_{\T^8}+\frac{1}{2}\omega_8)) =-\Upsilon_{\T^8}(\rho_{\T^8}-\frac{1}{2}\omega_8)$.
Hence, we have
\begin{align*}
\sum_{\beta\in\Xi}\Upsilon_{\T^8}(s_{\beta}(\rho_{\T^8}+\frac{1}{2}\omega_8))=0.
\end{align*}
This means that Jantzen's criterion is satisfied.

\,\,\, Q.E.D.

{\bf Remark} \,\,\,\, In the case of $k=7$, the non-existence of the homomorphism is proved in
[Boe-Collingwood 1990] for the regular integral case.

\subsection{${\bf E_8}$}

We fix an orthonormal basis $e_1,...,e_8$ in $\hhh^\ast$
such  that
\begin{multline*}
\Delta=\{\pm e_i\pm e_j\mid 1\leqslant i<j\leqslant
8\}\\
\cup \left\{\pm\frac{1}{2}\left(\sum_{i=1}^8(-1)^{n(i)}e_i\right)\left|
\mbox{$n(i)$ is either $0$ or $1$ for $3\leqslant i\leqslant 8$ and
 $\sum_{i=3}^8n(i)$ is odd.} \right.\right\}
\end{multline*}

We choose a positive system as follows.
\begin{multline*}
\Delta^+=\{ e_i\pm e_j\mid 1\leqslant i<j\leqslant
8\}\\
\cup \left\{\frac{1}{2}\left(e_1+\sum_{i=2}^8(-1)^{n(i)}e_i\right)\left|
\mbox{$n(i)$ is either $0$ or $1$ for $2\leqslant i\leqslant 8$ and
 $\sum_{i=2}^8n(i)$ is odd.} \right.\right\}
\end{multline*}

Put $\alpha_i=e_{i+1}-e_{i+2}$ for $1\leqslant i\leqslant 5$,
$\alpha_7=e_7+e_8$, and
$\alpha_7=\frac{1}{2}(e_1-e_2-e_3-e_4-e_5-e_6-e_7-e_8)$.
Then, $\Pi=\{\alpha_1,...,\alpha_8\}$ is the set of simple roots in
$\Delta^+$.
\[
 \begin{array}{ccccccccccccc} 1 & - & 2 & - & 3 & - & 4 & - & 5 & - & 7 & - & 8\\
& &  & & & & & & \mid & & & & \\
& & & &  & & & & 6 & & & & 
\end{array}
\]

We have:

\begin{thm} \mbox{}

(1) \,\,\,\,  Assume $k\in\{1,2,4,6,8\}$.
Then, $M_{\T^k}[-t]\subseteq M_{\T^k}[t]$ if and only if $t\in\nat$.

(2) \,\,\,\,  Assume $k\in\{3,5,7\}$.
Then, $M_{\T^k}[-t]\subseteq M_{\T^k}[t]$ if and only if $t\in\frac{1}{2}\nat$.
\end{thm}

\proof
For any $k$, we have $w_{\T^k}w_0=w_0w_{\T^k}$ since $w_0$ is contained
in the center of $W$.
Hence, Lemma 2.2.3 implies that $M_{\T^k}[-t]\subseteq M_{\T^k}[t]$,
only if $t\in\frac{1}{2}\nat$.

Next we consider the case of $t\in\nat$.
For $k\in\{1,2,3,4,6,8\}$, $\ppp_{\T^k}$ is even (cf.\ [Carter 1985]
p405-406).
In this case, Corollary 2.2.9 and Lemma 4.1.2 imply  $M_{\T^k}[-t]\subseteq M_{\T^k}[t]$ for  $t\in\nat$.
If $k=5$, we have  $d_5=\frac{9}{2}$.
So, $\rho_{\T^5}$ is not integral.
A basis of integral root system for $\rho_{\T^5}$ is
$\T^5\cup\{e_5+e_6\}$.
We see that the integral root system is of the type $A_1\times E_7$.
So, the problem is reduced to the case of $k=3$ in the type $E_7$.
Hence, $M_{\T^5}[-t]\subseteq M_{\T^5}[t]$ for $t\in\nat$ such that
$\rho_{\T^5}+t\omega_5$ is dominant and regular.
From Lemma 2.2.6, we have  $M_{\T^5}[-t]\subseteq M_{\T^5}[t]$ for all
$t\in\nat$.
The case of $k=7$ is similar to the case of $k=5$.
This time, a basis of the integral Weyl group of $\rho_{\T^7}$ is
$\T^7\cup\{\frac{1}{2}(e_1-e_2-e_3-e_4+e_5+e_6+e_7+e_8)\}$ and the
integral root system is of the type $A_1\times E_7$.
Hence, the problem is reduced to the case of $k=5$ in the type
$E_7$.

Next, we consider the case of $t\in\frac{1}{2}+\nat$.

First, we consider the case of $k=3$.
In this case, $\rho_{\T^3}+t\omega_3$ is not integral for $t\in\frac{1}{2}+\nat$.
A basis of integral root system for $\rho_{\T^3}+\frac{1}{2}\omega_3$ is
$\T^5\cup\{e_3+e_4\}$ and the
integral root system is of the type $D_8$. 
This time, the problem is reduced to the case of $k=3$ in the type
$D_8$.

Next, we consider the case of $k=5,7$.
In this case, $\rho_{\T^k}+t\omega_k$ is integral for
$t\in\frac{1}{2}+\nat$.
The special representation corresponding to a Richardson orbit can be
constructed as a MacDonald representation ([MacDonald 1972]).
If $k=5$ (resp.\ $k=7$), then the numder of positive roots in the Levi
part of $\ppp_{\T^k}$ is $14$ (resp.\ $22$).
Hence, the special representation occurs as component of  $S^{14}(\hhh^\ast)$ (resp.\
$S^{22}(\hhh^\ast)$) but not of $S^{d}(\hhh^\ast)$ for $d>14$ (resp.\
$d>22$).
Hence, from [Carter 1985] p417-418, we see the family of the special
representation consists of a single element.
Hence, Corollary 2.1.4, Lemma 2.2.6  and Lemma 2.2.7 inply that 
$M_{\T^k}[-t]\subseteq M_{\T^k}[t]$ for $t\in\frac{1}{2}+\nat$.

If $k=1,2,4,6,8$, we can prove that $M_{\T^k}[\frac{1}{2}]$ is irreducible via
Jantzen's criterion. So, we have $M_{\T^k}[-t]\not\subseteq M_{\T^k}[t]$ for $t\in\frac{1}{2}+\nat$.
We describe the computation  briefly.

First, we consider the case of $k=1$.
In this case $d_1=\frac{29}{2}$ and $\rho_{\T^1}+\frac{1}{2}\omega_1$ is integral.
We put $\gamma_1=e_1+e_2$ and  $\gamma_2=\frac{1}{2}(e_1+e_2+e_3-e_4+e_5-e_6+e_7-e_8)$.

Put $\Xi=\left\{\beta\in\left. \left(\Delta^{\T^1}\right)^+ \right|
 \langle\rho_{\T^1}+\frac{1}{2}\omega_1,\beta^\vee\rangle\in\nat-\{0\}\right\}$.
Then, $\gamma_1,\gamma_{2}\in\Xi$.
For $\beta\in\Xi-\{\gamma_1, \gamma_2\}$, we can find
 $\eta\in\Delta_{\T^1}$ such that $\langle\beta,\eta\rangle=0$.
Hence, we have
 $\Upsilon_{\T^1}(s_\beta(\rho_{\T^1}+\frac{1}{2}\omega_1))=0$.
Moreover, we have
$\Upsilon_{\T^1}(s_{\gamma_{1}}(\rho_{\T^1}+\frac{1}{2}\omega_1))
=\Upsilon_{\T^1}(\rho_{\T^1}-\frac{1}{2}\omega_1)$ and
$\Upsilon_{\T^1}(s_{\gamma_{2}}(\rho_{\T^1}+\frac{1}{2}\omega_1)) =-\Upsilon_{\T^1}(\rho_{\T^1}-\frac{1}{2}\omega_1)$.
Hence, we have
\begin{align*}
\sum_{\beta\in\Xi}\Upsilon_{\T^1}(s_{\beta}(\rho_{\T^1}+\frac{1}{2}\omega_1))=0.
\end{align*}
This means that Jantzen's criterion is satisfied.

Next,we consider the case of $k=2$.
In this case $d_2=\frac{19}{2}$ and $\rho_{\T^2}+\frac{1}{2}\omega_2$ is integral.
We remark that $\PP_{\T^k}^{++}\cap W(\rho_{\T^2}+\frac{1}{2}\omega_2)$
consists of the $3$ elements (say $\lambda_1,\lambda_2, \lambda_3$).
We put  $\lambda_1=\rho_{\T^2}+\frac{1}{2}\omega_2$ and $\lambda_3=\rho_{\T^2}-\frac{1}{2}\omega_2$.
The remaining element is $\lambda_2=4e_1-3e_2-5e_3+4e_4+3e_5+2e_6+e_7$.
We put $\gamma_1=e_1+e_2$, $\gamma_2=e_1+e_3$, $\gamma_3=e_1-e_6$, and $\gamma_4=\frac{1}{2}(e_1+e_2-e_3+e_4-e_5-e_6+e_7+e_8)$.
Put $\Xi=\left\{\beta\in\left. \left(\Delta^{\T^2}\right)^+ \right|
 \langle\rho_{\T^2}+\frac{1}{2}\omega_2,\beta^\vee\rangle\in\nat-\{0\}\right\}$.
Then, $\gamma_1,...,\gamma_{4}\in\Xi$.
For $\beta\in\Xi-\{\gamma_1,...,\gamma_{4}\}$, we can find
 $\eta\in\Delta_{\T^2}$ such that $\langle\beta,\eta\rangle=0$.
Hence, we have
 $\Upsilon_{\T^2}(s_\beta(\rho_{\T^2}+\frac{1}{2}\omega_2))=0$.
Moreover, we have
$\Upsilon_{\T^2}(s_{\gamma_{1}}(\rho_{\T^2}+\frac{1}{2}\omega_2)) =-\Upsilon_{\T^2}(\lambda_3)$,
$\Upsilon_{\T^2}(s_{\gamma_{2}}(\rho_{\T^2}+\frac{1}{2}\omega_2)) =\Upsilon_{\T^2}(\lambda_2)$,
$\Upsilon_{\T^2}(s_{\gamma_{3}}(\rho_{\T^2}+\frac{1}{2}\omega_2))
=\Upsilon_{\T^2}(\lambda_3)$, and
$\Upsilon_{\T^2}(s_{\gamma_{4}}(\rho_{\T^2}+\frac{1}{2}\omega_2))
 =-\Upsilon_{\T^2}(\lambda_2)$.
Hence, we have
\begin{align*}
\sum_{\beta\in\Xi}\Upsilon_{\T^2}(s_{\beta}(\rho_{\T^2}+\frac{1}{2}\omega_2))=0.
\end{align*}
This means that Jantzen's criterion is satisfied.

Next, we consider the case of  $k=4$.
In this case $d_4=\frac{11}{2}$ and $\rho_{\T^4}+\frac{1}{2}\omega_4$ is integral.
We remark that $\PP_{\T^k}^{++}\cap W(\rho_{\T^4}+\frac{1}{2}\omega_4)$
consists of the $7$ elements
$\lambda_1,...,\lambda_7$.
They are characterized as follows.
$\langle\lambda_1,\omega_4\rangle=10$,
$\langle\lambda_2,\omega_4\rangle=5$,
$\langle\lambda_3,\omega_4\rangle=2$,
$\langle\lambda_4,\omega_4\rangle=0$,
$\langle\lambda_5,\omega_4\rangle=-2$,
$\langle\lambda_6,\omega_4\rangle=-5$, and 
$\langle\lambda_7,\omega_4\rangle=-10$.
We have $\lambda_1=\rho_{\T^4}+\frac{1}{2}\omega_4$ and
$\lambda_7=\rho_{\T^4}-\frac{1}{2}\omega_4$.
We put $\gamma_1=e_1+e_2$, $\gamma_2=e_1+e_3$, $\gamma_3=e_1+e_4$,
$\gamma_4=e_1+e_5$, $\gamma_5=e_1+e_6$, $\gamma_6=e_1-e_5$,
$\gamma_7=e_1-e_7$, $\gamma_8=e_1-e_8$,
$\gamma_9=\frac{1}{2}(e_1+e_2+e_3-e_4-e_5+e_6+e_7-e_8)$, 
$\gamma_{10}=\frac{1}{2}(e_1+e_2-e_3-e_4-e_5+e_6+e_7+e_8)$,
$\gamma_{11}=\frac{1}{2}(e_1+e_2+e_3+e_4-e_5+e_6+e_7+e_8)$, and
$\gamma_{12}=\frac{1}{2}(e_1+e_2+e_3-e_4-e_5+e_6-e_7+e_8)$.

Put $\Xi=\left\{\beta\in\left. \left(\Delta^{\T^4}\right)^+ \right|
 \langle\rho_{\T^4}+\frac{1}{2}\omega_4,\beta^\vee\rangle\in\nat-\{0\}\right\}$.
Then, $\gamma_1,...,\gamma_{12}\in\Xi$.
For $\beta\in\Xi-\{\gamma_1,...,\gamma_{12}\}$, we can find
 $\eta\in\Delta_{\T^4}$ such that $\langle\beta,\eta\rangle=0$.
Hence, we have
 $\Upsilon_{\T^4}(s_\beta(\rho_{\T^4}+\frac{1}{2}\omega_4))=0$.
Moreover, we have
$\Upsilon_{\T^4}(s_{\gamma_{1}}(\rho_{\T^4}+\frac{1}{2}\omega_4)) =-\Upsilon_{\T^4}(\lambda_7)$,
$\Upsilon_{\T^4}(s_{\gamma_{2}}(\rho_{\T^4}+\frac{1}{2}\omega_4)) =\Upsilon_{\T^4}(\lambda_6)$,
$\Upsilon_{\T^4}(s_{\gamma_{3}}(\rho_{\T^4}+\frac{1}{2}\omega_4)) =-\Upsilon_{\T^4}(\lambda_4)$,
$\Upsilon_{\T^4}(s_{\gamma_{4}}(\rho_{\T^4}+\frac{1}{2}\omega_4))
 =\Upsilon_{\T^4}(\lambda_2)$,
$\Upsilon_{\T^4}(s_{\gamma_{5}}(\rho_{\T^4}+\frac{1}{2}\omega_4))
 =\Upsilon_{\T^4}(\lambda_7)$,
$\Upsilon_{\T^4}(s_{\gamma_{6}}(\rho_{\T^4}+\frac{1}{2}\omega_4))
 =-\Upsilon_{\T^4}(\lambda_6)$,
$\Upsilon_{\T^4}(s_{\gamma_{7}}(\rho_{\T^4}+\frac{1}{2}\omega_4))
 =\Upsilon_{\T^4}(\lambda_3)$,
$\Upsilon_{\T^4}(s_{\gamma_{8}}(\rho_{\T^4}+\frac{1}{2}\omega_4))
 =-\Upsilon_{\T^4}(\lambda_5)$,
$\Upsilon_{\T^4}(s_{\gamma_{9}}(\rho_{\T^4}+\frac{1}{2}\omega_4))
 =\Upsilon_{\T^4}(\lambda_4)$,
$\Upsilon_{\T^4}(s_{\gamma_{10}}(\rho_{\T^4}+\frac{1}{2}\omega_4))
 =-\Upsilon_{\T^4}(\lambda_2)$,
$\Upsilon_{\T^4}(s_{\gamma_{11}}(\rho_{\T^4}+\frac{1}{2}\omega_4))
 =\Upsilon_{\T^4}(\lambda_5)$, and
$\Upsilon_{\T^4}(s_{\gamma_{12}}(\rho_{\T^4}+\frac{1}{2}\omega_4))
 =-\Upsilon_{\T^4}(\lambda_3)$.
Hence, we have
\begin{align*}
\sum_{\beta\in\Xi}\Upsilon_{\T^4}(s_{\beta}(\rho_{\T^4}+\frac{1}{2}\omega_4))=0.
\end{align*}
This means that Jantzen's criterion is satisfied.

Next, we consider the case of $k=8$.
In this case $d_1=\frac{23}{2}$ and $\rho_{\T^1}+\frac{1}{2}\omega_8$ is integral.
We put $\gamma_1=e_1+e_7$ and  $\gamma_2=\frac{1}{2}(e_1+e_2+e_3-e_4-e_5+e_6+e_7-e_8)$.

Put $\Xi=\left\{\beta\in\left. \left(\Delta^{\T^8}\right)^+ \right|
 \langle\rho_{\T^8}+\frac{1}{2}\omega_8,\beta^\vee\rangle\in\nat-\{0\}\right\}$.
Then, $\gamma_1,\gamma_{2}\in\Xi$.
For $\beta\in\Xi-\{\gamma_1, \gamma_2\}$, we can find
 $\eta\in\Delta_{\T^8}$ such that $\langle\beta,\eta\rangle=0$.
Hence, we have
 $\Upsilon_{\T^8}(s_\beta(\rho_{\T^8}+\frac{1}{2}\omega_8))=0$.
Moreover, we have
$\Upsilon_{\T^8}(s_{\gamma_{1}}(\rho_{\T^8}+\frac{1}{2}\omega_8))
=\Upsilon_{\T^8}(\rho_{\T^8}-\frac{1}{2}\omega_8)$ and
$\Upsilon_{\T^8}(s_{\gamma_{2}}(\rho_{\T^8}+\frac{1}{2}\omega_8)) =-\Upsilon_{\T^8}(\rho_{\T^8}-\frac{1}{2}\omega_8)$.
Hence, we have
\begin{align*}
\sum_{\beta\in\Xi}\Upsilon_{\T^8}(s_{\beta}(\rho_{\T^8}+\frac{1}{2}\omega_8))=0.
\end{align*}
This means that Jantzen's criterion is satisfied.

Finally, we consider the case $k=6$.
In this case, we choose the following basis of the root system in order
to make computation easier.
$\alpha_1=e_7-e_8$, $\alpha_2=e_6-e_7$, $\alpha_3=e_5-e_6$,
$\alpha_4=e_4-e_5$, $\alpha_5=e_3-e_4$,
$\alpha_6=\frac{1}{2}(-e_1-e_2-e_3+e_4+e_5+e_6+e_7+e_8)$,
$\alpha_7=e_2-e_3$,
and $\alpha_8=e_1-e_2$.
In this case $d_6=\frac{17}{2}$ and $\rho_{\T^6}+\frac{1}{2}\omega_6$ is integral.
We put $\gamma_1=e_1+e_5$ and  $\gamma_2=\frac{1}{2}(e_1+e_2+e_3+e_4+e_5-e_6-e_7-e_8)$.
Put $\Xi=\left\{\beta\in\left. \left(\Delta^{\T^6}\right)^+ \right|
 \langle\rho_{\T^6}+\frac{1}{2}\omega_6,\beta^\vee\rangle\in\nat-\{0\}\right\}$.
Then, $\gamma_1,\gamma_{2}\in\Xi$.
For $\beta\in\Xi-\{\gamma_1, \gamma_2\}$, we can find
 $\eta\in\Delta_{\T^6}$ such that $\langle\beta,\eta\rangle=0$.
Hence, we have
 $\Upsilon_{\T^6}(s_\beta(\rho_{\T^6}+\frac{1}{2}\omega_6))=0$.
Moreover, we have
$\Upsilon_{\T^6}(s_{\gamma_{1}}(\rho_{\T^6}+\frac{1}{2}\omega_6))
=\Upsilon_{\T^6}(\rho_{\T^6}-\frac{1}{2}\omega_6)$ and
$\Upsilon_{\T^6}(s_{\gamma_{2}}(\rho_{\T^6}+\frac{1}{2}\omega_6)) =-\Upsilon_{\T^6}(\rho_{\T^6}-\frac{1}{2}\omega_6)$.
Hence, we have
\begin{align*}
\sum_{\beta\in\Xi}\Upsilon_{\T^6}(s_{\beta}(\rho_{\T^6}+\frac{1}{2}\omega_6))=0.
\end{align*}
This means that Jantzen's criterion is satisfied.
\,\,\, Q.E.D.

{\bf Remark} \,\,\,\, In the case of $k=1$, the non-existence of the homomorphism is proved in
[Boe-Collingwood 1990] for the regular integral case.


\setcounter{section}{5}
\setcounter{subsection}{0}

\section*{\S\,\, 5.\,\,\,\, Elementary homomorphisms}

Here, we explain that we can construct homomorphisms in the setting of
general parabolic subalgebras from the case of maximal parabolic
subalgebras. 

\subsection{A comparison result}

Here, we review some notion in [Matumoto 1993] \S 3.
Hereafter, $\gggg$ means a reductive Lie algebra over $\cpx$ and retain
the notations in \S 1.
We fix a subset $\T$ of $\Pi$.
For $\alpha\in\Delta$, we put
\begin{gather*}
 \Delta(\alpha)=\{\beta\in\Delta\mid\exists c\in\rel \,\,\,\, \beta|_{\aas}=c\alpha|_{\aas}\},\\
\Delta^+(\alpha)=\Delta(\alpha)\cap\Delta^+, \\
U_\alpha=\cpx S+\cpx\alpha\subseteq\hhd.
\end{gather*}

Then $(U_\alpha,\Delta(\alpha),\langle\,\,\,,\,\,\,\rangle)$ is a subroot system of $(\hhh^\ast, \Delta,\langle\,\,\,,\,\,\,\rangle)$.
The set of simple roots for $\Delta^+(\alpha)$ is denoted by $\Pi(\alpha)$.
 If $\alpha|_{\aas}\neq
0$, then there exists a unique $\tilde{\alpha}\in\Delta^+$ 
such that $\Pi(\alpha)=\T\cup\{\tilde{\alpha}\}$.
If $\alpha\in\Delta$ satisfies $\alpha|_{\aas}\neq 0$ and $\alpha=\tilde{\alpha}$, then we call $\alpha$ $\T$-reduced.
For $\alpha\in\Delta^+$, we denote by $W_\T(\alpha)$ the Weyl group of $(\hhh^\ast, \Delta(\alpha))$.  
Clearly, $W_\T\subseteq W_\T(\alpha)\subseteq W$.
We denote by $w^\alpha$ the longest element of $W_\T(\alpha)$.
We call $\alpha\in\Delta$ $\T$-acceptable if  $w^\alpha w_\T=w_\T w^\alpha$.
We denote by $\Delta_r^\T$ the set of  $\T$-reduced $\T$-acceptable roots.
Put $(\Delta^\T_r)^+=\Delta^+\cap\Delta^{\T}_r$.
For $\alpha\in \Delta^\T_r$, we define
\[\sigma_\alpha=w^\alpha w_\T=w_\T w^\alpha.\]
Clearly, ${\sigma_\alpha}^2=1$.
For $\alpha\in\Delta$, we put 
\[ V_\alpha=\{\lambda\in\ads\mid\langle\lambda,\alpha\rangle=0\}.\]

We denote by $\omega_\alpha\in\ads\subseteq\hhd$ the fundamental weight
for $\alpha$ with respect to the basis $\Pi(\alpha)=\T\cup\{\alpha\}$.
Namely $\omega_\alpha$ satisfies that $\langle\omega_\alpha,
\beta\rangle=0$ for $\beta\in\T$, $\langle \beta,\check{\alpha}\rangle=1$,
and $\omega_\alpha|_{\hhh\cap\ccc(\gggg(\alpha))}=0$.
Here, $\ccc(\gggg(\alpha))$ is the center of $\gggg(\alpha)$. 
We see that there is some positive real number $a$ such that
$\omega_\alpha=a\alpha|_{\aas}$, since
$\alpha|_{\hhh\cap\ccc(\gggg(\alpha))}=0 $.
Hence, we have $V_\alpha=\{\lambda\in\ads\mid \langle \lambda,\omega_\alpha\rangle=0\}$.

We can easily see:
\begin{lem}
Let $\alpha\in \Delta_r^\T$.  Then, we have

(1) \,\,\, $\sigma_\alpha$ preserves $\ads$.

(2) \,\,\, $\sigma_\alpha\in W(\T)$.  In particular, $\sigma_\alpha\rho_\T=\rho_\T$. 

(3) \,\,\,
$\sigma_{\alpha}\omega_\alpha=-\omega_\alpha$.

(4) \,\,\, $\sigma_\alpha|_{\aas}$ is the reflection with respect to $V_\alpha$.
\end{lem}

For $\alpha\in(\Delta^\T_r)^+$, we define
\begin{gather*}
\gggg(\alpha) =\hhh+\sum_{\beta\in\Delta(\alpha)}\gggg_\beta, \,\,\,\,\,\,\,\,\,
 \ppp_\T(\alpha) =\gggg(\alpha)\cap\ppp_{\T}.
\end{gather*}
Then, $\gggg(\alpha)$ is a reductive Lie subalgebra of $\gggg$ whose root system is $\Delta(\alpha)$ and  $\pps(\alpha)$ is a maximal parabolic subalgebra of  $\gggg(\alpha)$.

Put
$\rho(\alpha)=\frac{1}{2}\sum_{\beta\in\Delta^+(\alpha)}\beta$,
For $\nu\in\ads$, we denote by $\cpx_\nu$ the one-dimensional $U(\pps(\alpha))$-module corresponding to $\nu$.
For $\nu\in\ads$ we define a generalized Verma module for $\gggg(\alpha)$ as follows.
\begin{align*}
M_{\T}^{\gggg(\alpha)}(\rho_{\T}+\nu)=U(\gggg(\alpha))\otimes_{U(\pps(\alpha))}\cpx_{\nu-\rho(\alpha)}.
\end{align*}
Then, we have:
\begin{thm}
Let $\nu$ be an arbitrary element in $V_\alpha$ and let $c$ be either
 $1$ or $\frac{1}{2}$.
Assume that $M_{\T}^{\gggg(\alpha)}(\rho_{\T}-nc\omega_\alpha)\subseteq M_{\T}^{\gggg(\alpha)}(\rho_{\T}+nc\omega_\alpha)$ for all $n\in\nat$.
Then, we have
$M_{\T}(\rho_{\T}+\nu-nc\omega_\alpha)\subseteq M_{\T}(\rho_{\T}+\nu+nc\omega_\alpha)$ for all $n\in\nat$.
(We call the above homomorphism of
 $M_{\T}(\rho_{\T}+\nu-nc\omega_\alpha)$ into
 $M_{\T}(\rho_{\T}+\nu+nc\omega_\alpha)$ an elementary homomorphism. )
\end{thm}

\proof
Assume that $M_{\T}^{\gggg(\alpha)}(\rho_{\T}-nc\omega_\alpha)\subseteq
M_{\T}^{\gggg(\alpha)}(\rho_{\T}+nc\omega_\alpha)$ for all $n\in\nat$.
Remark that $\sigma_\alpha(\rho_{\T}+nc\omega_\alpha)=\rho_{\T}-nc\omega_\alpha$.
From Theorem 1.4.2, this implies that $\sigma_\alpha$ is a Duflo involution of the integral
Weyl group $W_{\rho_{\T}+nc\omega_\alpha}$ for a sufficiently large $n$.

Put $Q_{\alpha, n}=\{\nu\in V_\alpha\mid \Delta(\alpha)_{\rho_{\T}+nc\omega_\alpha}=\Delta_{\rho_{\T}+\nu+nc\omega_\alpha}\}$.
From Theorem 1.4.2, for sufficiently large $n\in\nat$ and $\nu\in Q_{\alpha, n}$, we
have $M_{\T}(\rho_{\T}+\nu-nc\omega_\alpha)\subseteq
M_{\T}(\rho_{\T}+\nu+nc\omega_\alpha)$, since we have $\sigma_\alpha(\rho_{\T}+\nu+nc\omega_\alpha)=\rho_{\T}+\nu-nc\omega_\alpha$.
We easily see $\nu-nc\omega_\alpha$ is strongly $\T$-antidominant for
all $n\in\nat$.
Applying Lemma 1.5.1 and the exactness of the translation functor, we
can remove the extra assumption that $n$ is sufficiently large.

On the other hand $Q_{\alpha,n}$ is Zarisky dense in $V_\alpha$ (Cf.\
[Matumoto 1993] Lemma 3.2.2 (1)).
Moreover, we can prove that 
$\{\nu\in \ads\mid M_{\T}(\rho_{\T}+\nu-\mu)\subseteq
M_{\T}(\rho_{\T}+\nu) \}$ is Zarisky closed in $\ads$ for each
$\mu\in\ads$ in the same way as [Lepowsky 1975b] Lemma 5.4.
Hence, for each $\nu\in V_\alpha$ and each $n\in\nat$, we
have  $M_{\T}(\rho_{\T}+\nu-nc\omega_\alpha)\subseteq
M_{\T}(\rho_{\T}+\nu+nc\omega_\alpha)$.
\,\,\,\, Q.E.D.

\medskip
\medskip

{\bf Remark} \,\,\,\, Taking this opportunity, I would like to fix an
error in [Matumoto 1993].
In page 269 line 18, the definitions of $\gggg(\alpha,c)$ and
$\ppp_S(\alpha,c)$ are incorrect.
$\gggg(\alpha,c)$ should be an abstract reductive Lie algebra associated
with the pair $(\hhh, \Delta_{\alpha,c})$.
$\ppp_S(\alpha,c)$ should be the standard parabolic subalgebra
corresponding to $\T$.
Since $\Delta_{\alpha,c}$ need not be closed under the addition in
$\Delta$, $\gggg(\alpha,c)$ need not be a subalgebra of $\gggg$.

\subsection{${\bf C_n}$ case}

As an example, we describe  elementary homomorphisms in the $C_n$ case.
Let $\gggg=\sss\ppp(n,\cpx)$.
We use the notation in the root system in 3.1.

Let $\kappa=(k_1,...,k_s)$ be a finite sequence of positive integers
such that 
$ k_1+\cdots+k_s\leqslant n$.
We put $k^\ast_i= k_1+\cdots k_i$ for $1\leqslant i\leqslant s$ and
$k^\ast_0=0$.
We define a subset $\T^\kappa$ of $\Pi$ as follows.
\[ \T^\kappa=\left\{\begin{array}{ll} 
\Pi-\{e_{k^\ast_i}-e_{k^\ast_i+1}| 1\leqslant i\leqslant s\} & \mbox{if $k^\ast_s<n$},\\
\Pi-(\{e_{k^\ast_i}-e_{k^\ast_i+1}| 1\leqslant i\leqslant s-1\}\cup \{2e_{n}\}) & \mbox{if $k^\ast_s=n$}
\end{array}\right. .\]
Then the corresponding standard Levi subalgebra $\llll_{\T^\kappa}$ is
isomorphic to
$\gl(k_1,\cpx)\oplus\gl(k_2,\cpx)\oplus\cdots\oplus\gl(k_s,\cpx)\oplus\sss\ppp(n-k^\ast_s,\cpx)$.
Here, we regard $\sss\ppp(0,\cpx)$ as a trivial Lie algebra $\{0\}$.
Obviously any proper subset of $\Pi$ is written as the form of $\T^\kappa$.

We put
$
 a_i=\sum_{j=1}^{k_i}e_{k^\ast_{i-1}+j} \,\,\,\,\, (1\leqslant
 i\leqslant s))$.
Then, $a_1,...,a_s$ form  a basis of $\aaa_{\T^\kappa}^\ast$.
We write $M_{\T^\kappa}[t_1,...,t_s]$ for
$M_{\T^\kappa}(\rho_{\T^\kappa}+t_1a_1+\cdots+t_sa_s)$ for $t_1,...,t_s\in\cpx$.

We easily have:

\begin{lem}\mbox{}

(1) \,\,\,\, If $k^\ast_s<n$, then 
\begin{align*}
(\Delta_r^{\T^\kappa})^+=\{e_{k^\ast_i}-e_{k^\ast_j+1}\mid 1\leq i\leqslant j< s,
 k_i=k_{j+1}\} \cup
\{e_{k^\ast_i}-e_{k^\ast_s+1}\mid 1\leq i\leqslant s\}
\end{align*}

(2) \,\,\,\, If $k^\ast_s=n$, then 
\begin{align*}
(\Delta_r^{\T^\kappa})^+=\{e_{k^\ast_i}-e_{k^\ast_j+1}\mid 1\leq i\leqslant j< s,
 k_i=k_{j+1}\} \cup
\{2e_{k^\ast_i}\mid 1\leq i\leqslant s\}
\end{align*}
\end{lem}

Combining [Boe 1985] 4.4 Theorem, Theorem 3.2.2, and Theorem 5.1.2, we
have
\begin{prop}
 
\,\, (1) \,\, Let $1\leqslant p< q\leqslant s$ be such that $k_i=k_j$.
If $t_p-t_q\in\nat$, we have
\begin{align*}
M_{\T^\kappa}\left(\rho_{\T^\kappa}+\sum_{\substack{1\leqslant i\leqslant s
 \\ i\neq p,q}} t_ia_i+t_qa_p+t_pa_q\right)\subseteq
 M_{\T^\kappa}\left(\rho_{\T^\kappa}+\sum_{1\leqslant i\leqslant s} t_ia_i\right).
\end{align*}

\,\, (2) \,\, Let $1\leqslant p\leqslant q$ be such that
 $3k_p>2(k_p+n-k^\ast_s)$.
If $t_p\in\nat$, we have
\begin{align*}
M_{\T^\kappa}\left(\rho_{\T^\kappa}+\sum_{\substack{1\leqslant i\leqslant s
 \\ i\neq p}} t_ia_i-t_pa_p\right)\subseteq
 M_{\T^\kappa}\left(\rho_{\T^\kappa}+\sum_{1\leqslant i\leqslant s} t_ia_i\right).
\end{align*}

\,\, (3) \,\, Let $1\leqslant p\leqslant q$ be such that
 $3k_p\leqslant 2(k_p+n-k^\ast_s)$ and $k_p$ is even.
If $t_p\in\frac{1}{2}\nat$, we have
\begin{align*}
M_{\T^\kappa}\left(\rho_{\T^\kappa}+\sum_{\substack{1\leqslant i\leqslant s
 \\ i\neq p}} t_ia_i-t_pa_p\right)\subseteq
 M_{\T^\kappa}\left(\rho_{\T^\kappa}+\sum_{1\leqslant i\leqslant s} t_ia_i\right).
\end{align*}

\end{prop}


\bigskip

\subsection*{References}

\mbox{}\hspace{5mm} 

{\sf [Barbasch-Vogan 1982]} D.\ Barbasch and D.\ A.\ Vogan Jr.\ , 
Primitive ideals and orbital integrals in complex classical groups,
{\it Math.\ Ann.\ } {\bf 259} (1982), 153--199.

{\sf [Barbasch-Vogan 1983]} D.\ Barbasch and D.\ A.\ Vogan Jr.\ ,
 Primitive ideals and orbital integrals in complex exceptional groups.
{\it J.\ Algebra} {\bf 80} (1983), 350--382.

{\sf [Baston 1990]} R.\ J.\ Baston, 
Verma modules and differential conformal invariants.
{\it J.\ Differential Geom.\ } {\bf 32} (1990), 851--898.

{\sf [Beilinson-Bernstein 1993]} A.\ Beilinson and J.\ Bernstein,
A proof of Jantzen conjectures,
I. M. Gelfand Seminar, 1-50,
{\it Adv.\ Soviet Math.\ }, {\bf 16} Part 1,
Amer.\ Math.\ Soc.\ , Providence, RI, 1993. 

{\sf [Bernstein-Gelfand-Gelfand 1971]} J.\ Bernstein, I.\ M.\ Gelfand, and S.\ I.\ Gelfand, Structure of representations generated by vectors of highest weight, {\it Funct.\ Anal. Appl.\ }{\bf 5} (1971), 1-8.

{\sf [Bernstein-Gelfand-Gelfand 1975]} J.\ Bernstein, I.\ M.\ Gelfand, and S.\ I.\ Gelfand, Differential operators on the base affine space and a study of $\gggg$-modules, in: ``Lie Groups and their Representations'', Wiley, New York, 1975, p 21-64.

{\sf [Bernstein-Gelfand 1980]} J.\ Bernstein and S.\ I.\ Gelfand, Tensor product of finite and infinite dimensional representations of semisimple Lie algebras, {\it Compos.\ Math.\ } {\bf 41} (1980), 245-285.

{\sf [Bien 1990]} F.\ Bien, ``$\diff$-Modules and Spherical Representations'' Math.\ Notes 39, Princeton University Press, Princeton, New Jersey, 1990.

{\sf [Boe 1985]} B.\ Boe, Homomorphism between generalized Verma modules, {\it Trans.\ Amer.\ Math.\ Soc.\ }{\bf 288} (1985), 791-799.

{\sf [Boe-Collingwood 1985]} B.\ Boe and D.\ H.\ Collingwood, A  comparison theory for the structure of induced representations, {\it J.\ of Algebra} {\bf 94} (1985), 511-545.

{\sf [Boe-Collingwood 1990]} B.\ Boe and D.\ H.\ Collingwood,
Multiplicity free categories of highest weight representations. I, II,
{\it Comm.\ Algebra} {\bf 18} (1990), 947--1032, 1033--1070.

{\sf [Boe-Enright-Shelton 1988]}
B.\ Boe, T.\ J.\ Enright, and B.\ Shelton,
Determination of the intertwining operators for holomorphically induced representations of Hermitian symmetric pairs.
{\it Pacific J.\ Math.\ } {\bf 131} (1988), 39--50.

{\sf [Borho-Jantzen 1977]} W.\ Borho and J.\ C.\ Jantzen, \"{U}ber primitive Ideale in der Einh\"{u}llenden einer halbeinfachen Lie-Algebra, {\it Invent.\ Math.\ }{\bf 39} (1977), 1-53.

{\sf [Borho-Kraft 1976]} W.\ Borho and H.\ Kraft, \"{U}ber die
Gelfand-Krillov-Dimension, {\it Math.\ Ann.\ }{\bf 220} (1976), 1-24.

{\sf [Brylinski-Kashiwara 1981]}J.-L.\ Brylinski and M.\ Kashiwara, 
Kazhdan-Lusztig conjecture and holonomic systems,
{\it Invent.\ Math.\ } {\bf 64} (1981), 387-410.

{\sf [Carter 1985]}
R.\ W.\ Carter, 
``Finite groups of Lie type.
Conjugacy classes and complex characters''. Pure and Applied Mathematics. A Wiley-Interscience Publication.
John Wiley \& Sons, Inc., New York, 1985. xii+544 pp. 

{\sf [Casian-Collingwood 1987]} 
The Kazhdan-Lusztig conjecture for generalized Verma modules.
{\it Math.\ Z.\ } {\bf 195} (1987), 581--600.

{\sf [Collingwood-Irving-Shelton 1988]} 
 D.\ H.\ Collingwood, R.\ S.\ Irving, and B.\ Shelton,
Filtrations on generalized Verma modules for Hermitian symmetric pairs.
{\it J.\ Reine Angew.\ Math.\ } {\bf 383} (1988), 54--86.

{\sf [Collingwood-McGovern 1993]} D.\ H.\ Collingwood and W.\ M.\ McGovern,
``Nilpotent orbits in semisimple Lie algebras''.
Van Nostrand Reinhold Mathematics Series.
Van Nostrand Reinhold Co., New York, 1993. xiv+186 pp. 

{\sf [Collingwood-Shelton 1990]} D.\ H.\ Collingwood and B.\ Shelton,
A duality theorem for extensions of induced highest weight modules,
{\it Pacific J.\ Math.\ } {\bf 146} (1990), 227--237.

{\sf [Dixmier 1977]} J.\ Dixmier, ``Enveloping Algebras'' North-Holland, Amsterdam and New York, 1977.

{\sf [Duflo 1977]} M.\ Duflo, Sur la classifications des id\'{e}aux primitifs dans l'alg\`{e}bre de Lie semi-simple, {\it Ann.\ of Math.\ } {\bf 105} (1977), 107-120.

{\sf [Gyoja 1994]} Highest weight modules and $b$-functions of semi-invariants,
{\it Publ.\ Res.\ Inst.\ Math.\ Sci.\ } {\bf 30} (1994), 353--400.

{\sf [Gyoja 2000]}  A.\ Gyoja, A duality theorem for homomorphisms between generalized Verma modules,
{\it J.\ Math.\ Kyoto Univ.\ }{\bf 40} (2000), 437--450.

{\sf [Hesselink 1978]}
W.\ H.\ Hesselink, 
Polarizations in the classical groups.
{\it Math.\ Z.\ } {\bf 160} (1978), 217--234.

{\sf [Jantzen 1977]}
J.\ C.\ Jantzen, 
Kontravariante Formen auf induzierten Darstellungen halbeinfacher Lie-Algebren,
{\it Math.\ Ann.\ } {\bf 226} (1977), 53--65.

{\sf [Joseph 1977]} A.\ Joseph, 
A characteristic variety for the primitive spectrum of a semisimple Lie
algebra, in :
Non-commutative harmonic analysis (Actes Colloq., Marseille-Luminy, 1976), pp. 102--118. Lecture Notes in Math., Vol. 587,
Springer, Berlin, 1977. 

{\sf [Joseph 1980a]} A.\ Joseph, 
Goldie rank in the enveloping algebra of a semisimple Lie algebra. I, 
{\it J.\ Algebra} {\bf 65} (1980), 269--283.

{\sf [Joseph 1980b]} A.\ Joseph, 
Goldie rank in the enveloping algebra of a semisimple Lie algebra. II,
{\it J.\ Algebra} {\bf 65} (1980), 284--306.

{\sf [Joseph 1983]} A.\ Joseph, On the classification of primitive ideals in the enveloping algebra of a semisimple Lie algebra, pp.\ 30-76 in: {\it Lecture Notes in Mathematics} No.\ 1024, Springer-Verlag, Berlin-Heidelberg-New York, 1983.

{\sf [Kazhdan-Lusztig 1979]}
D.\ Kazhdan and G.\ Lusztig, 
Representations of Coxeter groups and Hecke algebras.
{\it Invent.\ Math.\ } {\bf 53} (1979), 165--184.

{\sf [Knapp 2002]} A.\ W.\ Knapp, 
``Lie groups beyond an introduction''
Second edition. Progress in Mathematics, 140.
Birkh\"{a}user Boston, Inc., Boston, MA, 2002. xviii+812 pp.

{\sf [Lepowsky 1975a]} J.\ Lepowsky,
Conical vectors in induced modules,
{\it Trans.\ Amer.\ Math.\ Soc.\  } {\bf 208} (1975), 219-272.

{\sf [Lepowsky 1975b]} J.\ Lepowsky, Existence of conical vectors in induced modules, {\it Ann.\ of Math.\ }{\bf 102} (1975), 17-40.

{\sf [Lepowsky 1976]} J.\ Lepowsky,
Uniqueness of embeddings of certain induced modules,
{\it Proc.\ Amer.\ Math.\ Soc.\ } {\bf 56} (1976), 55--58.

{\sf [Lepowsky 1977]} J.\ Lepowsky, Generalized Verma modules, the Cartan-Helgason theorem, and the Harish-Chandra homomorphism, {\it J.\ Algebra} {\bf 49} (1977), 470-495.

{\sf [MacDonald 1972]} I.\ G.\ MacDonald, Some irreducible representations of Weyl groups, {\it Bull.\ London Math.\ Soc.\ } {\sf 4} (1972), 148-150.

{\sf [Matumoto 1993]} H.\ Matumoto, On the existence of homomorphisms
between scalar generalized Verma modules, in: Contemporary
Mathematics, 145, 259-274, Amer. Math. Soc., Providence, RI,
1993. 

{\sf [Matumoto 2003]} H.\ Matumoto,
The edge-of-wedge type embeddings of dereived functor modules for the
type A classical groups, preprint math.RT/0305329

{\sf [McGovern 1994]}  W.\ M.\ McGovern,
A remark on differential operator algebras and an equivalence of categories,
{\it Compo.\ Math.\ }{\bf 90} (1994), 305--313.

{\sf [Soergel 1987]} W.\ Soergel, Universelle versus relative Einh\"{u}llende:  \,\,\,\, Eine geometrische Untersuchung von Quotientienten von universellen Einh\"{u}llenden halbeinfacher Lie-Algebren, Dissertation, der Universit\"{a}t Hamberg, 1987.

{\sf [Verma 1968]} D.\ N.\ Verma, Structure of certain induced representations of complex semisimple lie algebras, {\it Bull.\ Amer.\ Math.\ Soc.\ }{\bf 74} (1968), 160-166.

{\sf [Vogan 1978]} D.\ A.\ Vogan Jr., Gelfand-Kirillov dimension for Harish-Chandra modules, {\it Invent.\ Math.\ } {\bf 48} (1978), 75-98.

{\sf [Vogan 1984]} D.\ A.\ Vogan Jr., Unitarizability of certain series of representations, {\it Ann.\ of Math.\ } {\bf 120} (1984), 141-187.

{\sf [Vogan 1986]} D.\ A.\ Vogan Jr., The orbit method and primitive
ideals for semisimple Lie algebras, in :
Lie algebras and related topics (Windsor, Ont., 1984), p281-316,
CMS Conf. Proc., 5,
Amer. Math. Soc., Providence, RI, 1986.

{\sf [Vogan 1988]} D.\ A.\ Vogan Jr., Irreducibilities of discrete series representations for semisimple symmetric spaces, {\it Adv.\  Stud.\ in Pure Math.\ } vol.\ 14, Kinokuniya Book Store, 1988, 381-417.

{\sf [Vogan 1990]} D.\ A.\ Vogan Jr.,
Dixmier algebras, sheets, and representation theory, in:
Operator algebras, unitary representations, enveloping algebras, and invariant theory (Paris, 1989), 333--395, Progr. Math., 92,
Birkh\"{a}user Boston, Boston, MA, 1990. 

{\sf [Yamashita 1986]} H.\ Yamashita, On Whittaker vectors for generalized Gelfand-Graev representations of semisimple Lie groups, {\it J.\ Math.\ Kyoto Univ.\ } {\bf 26} (1986), 263-298.

\end{document}